\documentclass[12pt]{article}
\usepackage[all]{xy}

\usepackage{epic}
\usepackage{amsmath}[1996/11/01]
\usepackage{amssymb,amsthm,amsfonts,latexsym,epsfig,graphics}
 \def\beql#1#2\eeql{\begin{equation}\label{#1}#2\end{equation}}

\advance \textheight by -1 \topmargin
\advance \textheight by -1 \headheight
\advance \textheight by -1 \headsep
\advance \textheight by -1 \footskip
\vsize = \textheight
\hsize = \textwidth
\DeclareMathOperator{\Fix}{Fix}
\DeclareMathOperator{\Min}{Min}

\DeclareMathOperator{\Stab}{Stab}

\DeclareMathOperator{\Hom}{Hom}
\DeclareMathOperator{\diag}{diag}
\DeclareMathOperator{\End}{End}
\DeclareMathOperator{\GL}{GL}
\DeclareMathOperator{\Aut}{Aut}
\DeclareMathOperator{\Alt}{Alt}

\DeclareMathOperator{\im}{im}

\DeclareMathOperator{\SL}{SL}
\DeclareMathOperator{\PSL}{PSL}
\DeclareMathOperator{\Cl}{Cl}

\newtheorem{theorem}{Theorem}[section]

\newtheorem{prop}[theorem]{Proposition}

\newtheorem{cor}[theorem]{Corollary}

\newcommand{\bew}{\noindent\underline{Proof.}\ }
\newtheorem{rem}[theorem]{Remark}
\newtheorem{lemma}[theorem]{Lemma}

\newtheorem{definition}[theorem]{Definition}
\newtheorem{defn}[theorem]{Definition}

\newcommand{\Z}{{\mathbb{Z}}}
\newcommand{\Q}{{\mathbb{Q}}}
\newcommand{\F}{{\mathbb{F}}}

\newcommand{\R}{{\mathbb{R}}}

\newcommand{\trace}{\mbox{trace}}

\newcommand{\eb}{\phantom{zzz}\hfill{$\square $}\smallskip}

\begin{document}
\begin{center}
{\Large {\bf On automorphisms of extremal even unimodular lattices of dimension 48.}}\\
\vspace{1.5\baselineskip}
{\em Gabriele Nebe} \\
\vspace*{1\baselineskip}
Lehrstuhl D f\"ur Mathematik, RWTH Aachen University\\
52056 Aachen, Germany \\
 nebe@math.rwth-aachen.de \\
\vspace{1.5\baselineskip}
\end{center}

{\sc Abstract.}
{\small The automorphism groups of the three known extremal even unimodular 
lattices of dimension 48 and the one of dimension 72  are computed using the
classification of finite simple groups. Restrictions on the possible automorphisms
of 48-dimensional extremal lattices are obtained. 
We classify all extremal lattices of dimension 48 having an automorphism of order $m$ with 
$\varphi(m) > 24$. In particular the lattice $P_{48n}$ is the unique extremal 48-dimensional
lattice that arises as an ideal lattice over a cyclotomic number field.
}

\section{Introduction}

Let $L$ be an even unimodular lattice in Euclidean $n$-space $(\R^n,(,))$, so 
$(x,x) \in 2\Z $ for all $x\in L$ and $L=L^{\#} =\{ x\in \R^n \mid (x,L) \subseteq \Z \} $.
Then the  theory of modular forms allows to upper bound the minimum 
$$\min (L) := \min \{ (x,x) \mid 0\neq x\in L \}
 \leq 2+ 2 \lfloor \frac{n}{24} \rfloor. $$ 
Extremal lattices are those even unimodular lattices $L$ that
achieve equality.
Of particular interest are
extremal even unimodular lattices $L$ in the jump dimensions, the multiples of 24.
There are only five extremal lattices known in jump dimensions: 
The Leech lattice $\Lambda _{24}$ of dimension 24, three lattices $P_{48p}$,
$P_{48q}$, $P_{48n}$ of dimension 48 and one lattice $\Gamma _{72}$ of dimension 72
(\cite{SPLAG}, \cite{cyclo}, \cite{dim72}). 
These five lattices realise the maximal known sphere packing density in these dimensions.
\\ 
Up to dimension 24 one knows
all even unimodular lattices, in particular 
the Leech lattice is the unique extremal lattice in dimension 24. 
A complete classification of all even unimodular lattices in dimension 48 seems to 
be impossible.
In the present paper we try to narrow down the possible automorphisms
of extremal 48-dimensional lattices $L$.
It turns out that all primes that can occur as an order of some automorphism
already occur for one of the three known 
examples. These primes are $47$, $23$ and all primes $\leq 13$. 
Explicit computations allow to find all extremal lattices $L$
 with an automorphism $\sigma $ of order $m$
such that $\varphi (m) > 24$ as well as those, where $\sigma$ has order 46 
and $\sigma ^{23} \neq -1$. 
The first section uses the classification of finite simple groups to 
prove the structure of the automorphism groups of the three known 48-dimensional
extremal lattices, as well as for $\Gamma _{72}$. 
We then continue to find all possible primes $p$ that might occur as the order of some 
automorphism of $L$ together with their fixed lattice $F$.
It turns out that $\dim (F) \leq 22$ for odd primes  $p$ and $F = 0, \sqrt{2}\Lambda _{24},
$ or $\sqrt{2} O_{24}$ for $p=2$. 
This allows to conclude that the minimal polynomial of any automorphism $\sigma \in \Aut(L)$
of order $m$ is divisible by the $m$-th cyclotomic polynomial $\Phi _m$ of degree
$\varphi(m)$. Using ideal lattices as introduced in \cite{ideal} this allows us to 
apply extensive number theoretic computations in Magma to find all lattices $L$ 
that have some automorphism $\sigma \in \Aut(L)$ of order $m$ with $\varphi (m) > 24$.
More precisely we prove the following theorem.

\begin{theorem}
Let $L$ be an extremal even unimodular lattice of dimension $48$
and $\sigma \in \Aut (L) $ of order $m$ such that $\varphi (m) > 24$ 
Then one of 
\begin{itemize}
\item $m=120$ and $L\cong P_{48n}$
\item $m=132$ and $L\cong P_{48p}$
\item $m=69$ and $L\cong P_{48p}$ 
\item $m=47$  and $L\cong P_{48q}$ 
\item $m=65$  and $L\cong P_{48n}$
\item $m=104$  and $L\cong P_{48n}$
\end{itemize}
\end{theorem}

\section{Bounds on the Hermite function.} \label{Bounds}

This section recalls some basic notions in the geometric theory of lattices. 
For more details the reader is referred to the textbook \cite{Martinet}.
The main purpose is to state the table displaying the bounds on the Hermite function
obtained from \cite{Elkies}.

Let $L = \bigoplus _{i=1}^n \Z B_i$ be a lattice in Euclidean space $(\R^n,(,))$
for some basis $B:=(B_1,\ldots , B_n)$ with 
Gram matrix $G(B):=((B_i,B_j)) _{i,j=1}^n $.
The determinant of $L$ is the determinant of
any of its Gram matrices, $\det (L) := \det (G(B)) $.
We denote by 
$$\Min (L):= \{ \ell \in L \mid (\ell,\ell) = \min (L ) \} .$$
the set of minimal vectors of $L$.
Its cardinality is known as the kissing number, as this is the number of 
spheres in the lattices sphere packing that touch one fixed sphere. 
The density of this sphere packing is maximal, if 
$\gamma (L)$ is maximal, where
the Hermite function  $\gamma $ on the space of similarity classes of 
$n$-dimensional lattices associates to $L$ the value
$$\gamma (L) := \frac{\min(L)}{\det (L) ^{1/n} } $$
The 
{\bf Hermite constant}  is
$$\gamma _n := \max \{ \gamma (L) \mid L \mbox{ is an } n- \mbox{dimensional lattice} \} .$$
Explicit values for $\gamma _n$ are only known for $n\leq 8$ and $n=24$. 
The best known upper bounds on $\gamma _n$ are given in \cite{Elkies}. 
Note that Cohn and Elkies work with the center density 
$\delta _n$. One gets $\gamma _n = 4 \delta _n^{2/n} $. 
The table below gives upper bounds $b_n$ for $\gamma _n$:
$$\begin{array}{|l|l||l|l||l|l||l|l||l|l|}
\hline
n & b_n & n & b_n & n & b_n & n & b_n & n & b_n \\
\hline
7 & 1.8115 & 13 & 2.6494 & 19 & 3.3975 & 25 & 4.1275 & 31 & 4.8484 \\
8 & 2 & 14 & 2.7759 & 20 & 3.5201 & 26 & 4.2481 & 32 & 4.9681 \\
9 & 2.1327 & 15 & 2.9015 & 21 & 3.6423 & 27 & 4.3685 & 33 & 5.0877 \\
10 & 2.2637 & 16 & 3.0264 & 22 & 3.7641 & 28 & 4.4887 & 34 & 5.2072 \\
11 & 2.3934 & 17 & 3.1507 & 23 & 3.8855 & 29 & 4.6087 & 35 & 5.3267 \\
12 & 2.5218 & 18 & 3.2744 & 24 & 4.0067 & 30 & 4.7286 & 36 & 5.4462 \\
\hline
\end{array} $$
The {\bf automorphism group} $$\Aut (L) := \{ \sigma \in O(\R^n,(,)) \mid \sigma (L) = L \} $$
respects the length and hence acts on $\Min (L)$. Taking matrices with respect to
the lattice basis $B$, we obtain $\Aut(L) \leq \GL_n(\Z )$ is a finite integral matrix group.

\section{Extremal even unimodular lattices of dimension 48} 

\subsection{Ternary codes and unimodular lattices of dimension 48}

Two of the known extremal even unimodular lattices of dimension 48 are 
have a canonical construction as 2-neighbors of code lattices of extremal ternary codes. 
Let me recall Sloane's construction  A: 

\begin{definition}(see for instance \cite[Chapter 7]{SPLAG})
Let $(e_1,\ldots , e_n)$ be a $p$-frame in $\R^n$, i.e. $(e_i,e_j) = p \delta _{ij} $.
For any code $C\leq (\Z/p\Z)^n $ the code lattice is
$A_p (C) := \{ \frac{1}{p} \sum c_i e_i \mid 
(\overline{c}_1, \ldots , \overline{c} _n) \in C \} $
where $\overline{x} := x+p\Z \in \Z/ p\Z $. 
\end{definition} 

Note that $A_p(C) $  contains the vectors $e_i$ of norm $p$.
To increase the minimum of the lattice one usually passes to 
a neighbor lattice as follows:  

\begin{definition} (\cite{Kneser})
Let $L$ be an integral lattice and $v\in L\setminus 2L^{\#} $ such that $(v,v) $ is a multiple of $4$. 
Then 
$$L^{(v),2} := \langle \{ \ell \in L\mid (v,\ell ) \mbox{ even } \} \cup \{ \frac{v}{2} \} \rangle _{\Z }$$
is called the {\em 2-neighbor} of $L$  defined by $v$. 
\end{definition}

Note that the 2-neighbor is an integral lattice with the same determinant as $L$.
If $L$ is unimodular, then any unimodular lattice $N$ such that $[L:N\cap L] =2$ is obtained as some 2-neighbor of
$L$.

Recall that the minimum weight of a self-dual ternary code $C$ of length 48 cannot exceed 15.
$C$ is called {\em extremal}, if $d(C) = 15$ (see for instance \cite{Koch}). 
Extremal codes always contain a vector of weight 48, so we may replace $C$ by some
equivalent code to obtain that the all-ones vector ${\bf 1} = (1,\ldots , 1) \in C$.

\begin{theorem}\cite{Koch} \label{Koch}
Let $C$ be an extremal self-dual ternary code of length $48$ containing the all-ones vector.
Then $\Lambda( C):=A_3(C)^{(v),2} $
is an extremal even unimodular lattice, where $v= \frac{1}{3} (e_1+\ldots + e_{48})\in A_3(C)$.
\end{theorem}

The extremal ternary self-dual codes of length 48 are not completely classified yet. 
One knows two equivalence classes of such codes
(see \cite{SPLAG}): the Pless code $P_{48}$ and the 
extended quadratic residue code $Q_{48}$. In fact \cite{LEN48} shows that these 
are the only such codes that have an automorphism of prime order $p\geq 5$:

\begin{theorem} \cite{LEN48}
Let $C$ be an extremal self-dual ternary code of length $48$ 
 such that $|\Aut(C)| $ is divisible
by some prime $p\geq 5$. Then $C\cong Q_{48} $ or $C\cong P_{48}$.
We have $\Aut(Q_{48}) \cong \SL_2(47)$ and  $\Aut (P_{48} ) \cong (\SL_2(23) \times C_2) :2 $.
\end{theorem}

The last theorem gives a characterization of the extremal lattices obtained 
from extremal ternary codes:

\begin{theorem} 
An extremal even unimodular lattice $L$ of dimension $48$ is of the form $L=\Lambda (C)$ 
for some extremal ternary code $C$ if and only if  
 there is some $\beta \in L$ with $(\beta ,\beta ) = 12$ such that 
$N_6(\beta ) := \{ x\in \Min (L) \mid (x,\beta ) = 6 \}$ has cardinality $94$. 
\end{theorem}

\bew
We first remark the well known fact that any $2$-neighbor $M$ of $L$ has minimum  $\geq 3$
and that the pairs of vectors $\pm v$ of norm 3 in $M$ are pairwise orthogonal:
This follows because any two vectors $v,w\in M \setminus L$ satisfy that both 
vectors $v\pm w\in M\cap L $ and hence 
$$(v \pm w ,v\pm w)  = (v,v) + (w,w) \pm 2 (v,w ) \geq 6 \mbox{ or } v \pm w = 0 .$$
Assume that there is such a $\beta \in L$. Then the neighbor $M:=L^{(\beta ),2} $ 
contains a 3-frame 
$$ \Min (M) = \{ \pm \frac{\beta }{2} \} \cup \{ x -\frac{\beta }{2} \mid x\in N_6(\beta ) \} $$
and so $M$ is a code lattice, $M=A_3(C)$ for some $C=C^{\perp } \leq \F_3^{48}$. 
Clearly $M$ is an odd unimodular lattice with even sublattice $M\cap L$. 
Since $\min (M\cap L) \geq \min (L)  = 6$, the code $C$ contains no words of weight
$\leq 12$, so $C$ is extremal.
\\
For $L=\Lambda (C) $ the vector $\beta = 2 e_1 \in L$ satisfies $(\beta , \beta )=12$ and
$N_6(\beta ) = \{ e_1 \pm e_j \mid j=2,\ldots ,48 \} $.
\eb

\subsection{The automorphism groups of the three known lattices}\label{aut48}

In the literature one finds three extremal even unimodular lattices of dimension 48. 
Two of them are constructed from the two known extremal ternary codes of 
length 48 as described in \cite{Koch} (see Theorem \ref{Koch}). 
They appear in \cite[Chapter 7, Example 9]{SPLAG}, where the authors refer to 
personal communications with J. Thompson for the structure of the automorphism group 
of these lattices. Since the description there is slightly incorrect and no explicit construction
of the automorphism groups is available in the literature, the construction of the
automorphism groups is given in the next theorem.
A third lattice, $P_{48n}$ has been found by the author in \cite{cyclo}, where it was 
proved that the normalizer of the subgroup $\SL_2(13)$ in $\Aut (P_{48n})$ is
$(\SL _2(13) {\sf Y} \SL_2(5) ).2^2 $. 
Using the classification of finite simple groups  one may obtain the full automorphism group
of these three lattices:

\begin{theorem}
\begin{tabular}{ccccc}
$\Aut(P_{48p}) \cong $ & $(\SL _2(23) \times S_3) : 2$ & of order & $72864 = 2^53^211\cdot 23  $ \\
 $\Aut (P_{48q}) \cong $ & $  \SL _2(47) $ & of order & $ 103776 =2^53\cdot23\cdot 47 $ \\ 
$\Aut (P_{48n}) \cong $ & $ (\SL _2(13) {\sf Y} \SL_2(5) ).2^2 $ & of order & $ 524160 = 2^73^25\cdot 7\cdot 13 $ 
\end{tabular}
\end{theorem}

\bew
Let $L$ be one of the three 48-dimensional extremal even unimodular lattices from the theorem
and let  $G:=\Aut (L)$ be its automorphism group. 
Then by construction $G$ contains the corresponding group $U$ from above as a subgroup.
Explicit matrices generating these subgroups can be 
obtained from the database of lattices \cite{LAT}.
Let $q =23,47,13$ be the largest prime divisor of $|U|$.  Then $U$ contains a normal subgroup
$\SL_2(q)$.
\\
{\bf (A)}
We first show that $U= N_G(\SL _2(q)) $:
\\
For $L=P_{48n}$, $q=13$ this is shown in \cite[Theorem 5.3]{cyclo}. 
\\
So let $L=P_{48q}$, $q=47$. Then $N:=\SL_2(47)$ acts on $L$ with endomorphism ring $\Z[\frac{1+\sqrt{-47}}{2} ]$. 
One computes that $L$ has a unique $\Z N$-sublattice $X$ of index 2. Since $N$ is perfect it fixes
all sublattices between $X^{\#}$ and $X$ and we obtain 3 unimodular lattices 
$$X < L,L',M < X^{\#} $$
with $L'$ even of minimum $4$ and $M$ odd with $\min (M) = 3$. 
The minimal vectors of $M$ form a 3-frame and hence $M= A_3(C) $ for some extremal ternary self-dual code $C$.
With MAGMA we compute $\Aut (C) \cong \SL_2(47)$ (in fact $C\cong Q_{48}$).
Since $\Aut(M)$ permutes the $48$ pairs of minimal vectors of $M$ we obtain 
$$\Aut (M) = \Stab _{C_2\wr S_{48}} (M) = \Aut(C) \cong \SL _2(47) .$$
The normalizer in $G=\Aut (L)$ of $N$ has to act on the unique sublattice $X$ and hence
permutes the lattices $\{ L,L',M \}$. These are pairwise non-isometric, so $N_G(N)$ 
stabilizes all three lattices, so
$N_G(N) = \Aut(M) = N \cong \SL_2(47) $. 
\\
Now let $L=P_{48p}$. Let $N:=\SL_2(23) \leq \Aut(L) $. Then $L$ has three $\Z N$-sublattices $X_i$ of index 2. 
The group $U$ permutes these three lattices $X_i$ transitively, so it is enough to show that 
$\Stab _{N_G(N)} (X_1) \leq U$. As before there are three unimodular lattices $X_1 \leq L,L',M \leq X_1^{\#}$ 
between $X_1$ and its dual lattice. $L = P_{48p}$ is even of minimum 6, $L'$ is even of minimum 4, and
$M$ is an odd lattice containing a 3-frame of minimal vectors. 
Again $M=A_3(C)$ for some extremal ternary self-dual code $C$.
With MAGMA we compute $\Aut (C) \cong (\SL_2(23) \times C_2):2$ (in fact $C\cong P_{48}$).
This group is isomorphic to $\Aut (M)$, fixes the even sublattice $X_1$ of $M$ and hence is 
isomorphic to  $\Stab _{N_G(N)} (X_1)$. In fact this is the easiest way to construct 
$N_G(N)$ and $\Aut(L)$. 
\\
{\bf (B)} 
Now  assume that $G=\Aut(L)$ is imprimitive: \\
First assume that $U\cong \SL_2(47)$ contained in $G = \Aut(L)$, so $L=P_{48q}$.
The group $U$ has a unique maximal subgroup of index $\leq 48$, this is $C_{47}:C_{46} $ and
of index 48 in $U$. The $U$-invariant $3$-frame constructed above gives rise to
the corresponding monomial representation of $U$. 
 So if $G=\Aut(L)$ is imprimitive,
then $U\leq G \leq C_2\wr S_{48} = \Aut (F)$, where $F \cong \sqrt{3} \Z^{48} $ is the lattice 
generated by the $U$-invariant 3-frame. 
Let $F_0$ be the even sublattice of $F$. Then $L$ contains $F_0$ of index $3^{24} 2$. 
Since $G$ stabilizes $L$ and $F$, it also stabilizes the even sublattice $X=A_3(C)_0 $ of the 
code lattice constructed above, 
$$ F_0 \underbrace{\leq }_{3^{24}} A_3(C)_0 \underbrace{\leq}_{2} L .$$
The same argument as above shows that this implies that $G\leq \Aut(C) \cong U$.
\\
Now let $U=(\SL_2(23) \times S_3):2 \leq G = \Aut (L)$ with $L=P_{48p}$. 
Since the derived subgroup $\SL_2(23) \times C_3$ of $U$ is rational irreducible, the 
representation is not induced from a 24-dimensional representation of a (normal) subgroup of index 2.
The composition factors for the subgroup $\SL_2(23)$  are of rational dimension 24, so if 
the representation of $U$ is rational imprimitive, then it is induced from a subgroup of
$U$ not containing $\SL_2(23)$.
The unique such subgroup of index $\leq 24$ has index 24 and is the normalizer $N$ in $U$ of a
23-Sylow subgroup of $U$, $N\cong (C_{46} : C_{11} \times S_3).2$. The restriction of the natural
representation of $U$ to $N$ has no composition factor of degree 2, so finally this shows that 
$U$ is a primitive subgroup of $\GL_{48}(\Q )$.
\\
The group $\SL_2(13)$ has no maximal subgroup of index dividing 48. 
From this one easily concludes that 
the group
$U = (\SL _2(13) {\sf Y} \SL_2(5) ).2^2$ is a 
primitive rational matrix group and so is $\Aut (P_{48n})$.
\\
{\bf (C)} 
Now assume that $G = \Aut(L)$ is a primitive subgroup of $\GL_{48}(\Q) $: \\
Then  all abelian normal subgroups of $G$ are cyclic. 
A theorem of Ph. Hall (see \cite[p. 357]{Hup}) 
classifies all $p$-groups whose abelian characteristic subgroups are cyclic
(these are central products of extraspecial $p$-groups with cyclic, dihedral, quaternion or quasidihedral groups). 
The relevant 2-groups (resp. 3-groups) that have a rational representation of dimension dividing 16 (resp. 6)
do not allow a non-trivial $\SL_2(q)$ action by automorphisms ($q=13,23,47$). 
Therefore we know that $O_2(G) = \langle - I_{48} \rangle $ and $O_3(G) = O_3(U) $($= C_3 $ or $1$) 
centralise $\SL_2(q)$.
If $G\neq U$ then $G$ has to contain a quasisimple proper overgroup of $\SL_2(q)$. 
The tables in Hiss and Malle \cite{HM} (which use the classification of finite simple groups)
exclude such overgroups inside $\GL_{48}(\Q )$. 
\eb

\begin{rem} 
The rational normaliser of the group $G= (\SL_2(13) {\sf Y} \SL_2(5)) .2^2 = \Aut(P_{48n})$  acts 
transitively on the $G$-invariant lattices in $\Q P_{48n}$ (see \cite{norm}).
So all $\Z G$-lattices in this space are similar to $P_{48n}$, therefore $G$ is maximal finite in $\GL_{48}(\Q )$.
\end{rem}

\subsection{Identifying new lattices} 

It is very hard, if not impossible with the present computing power, to check isometry of two 
extremal 48-dimensional even unimodular lattices, if nothing else but the Gram matrix is given. 
In this subsection I describe the computation of an explicit isometry between such lattices 
where I used the knowledge of a certain subgroup $U$ of order 48 of the automorphism group. 

In \cite{Miezaki} the authors construct two extremal even unimodular lattices 
$A_{14}(C_{14,48}) $ and $A_{46}(C_{46,48})$ of dimension 48
as code lattices of self-dual codes of length 48.
From the construction one also obtains a certain subgroup $U $ of order 48
(the monomial automorphism group of the codes $C_{2p,48} \leq \Z/(2p\Z) ^{48} $, 
with monomial entries $\pm 1$) which can be read off from the construction of the
codes with weighing matrices. With the help of this group $U$ we can construct an
explicit isometry $A_{14}(C_{14,48}) \cong P_{48n} $, $A_{46}(C_{46,48}) \cong P_{48p}$.

\begin{theorem}
 $A_{14}(C_{14,48}) \cong P_{48n} $ and $A_{46}(C_{46,48}) \cong P_{48p}$.
\end{theorem}

\bew
Let $L$ be one of the code lattices $A_{14}(C_{14,48})$ or $A_{46}(C_{46,48})$ and 
$U$ be the known subgroup of $\Aut(L)$ coming from the construction with codes. 
The group $U$ determines a sublattice $M \leq L$ as the full 
preimage of the fixed space $\Fix _U (L/2L) $ of the action of $U$ on $L/2L$.
The lattice $M$ contains $2L$ of index $2^2$, it has minimum norm 16 and kissing number 7200 and it 
is spanned by its minimal vectors. The automorphism group of $M$ has order 96. 
\\
The group $U$ has a normal subgroup $\langle \sigma \rangle $ of order 3. 
For both lattices $\Fix _L(\sigma ) \cong \sqrt{3} (E_8\perp E_8)$ and the orthogonal lattice
is a 32-dimensional lattice of determinant $3^{16}$ and minimum norm 6. 
Its automorphism group has order $2^83^9$ for $A_{14}(C_{14,48}) $ 
and $2^83^65$ for  $A_{46}(C_{46,48})$. 
\\
To find candidates among the three known lattices that might be isometric to $L$ we 
first find suitable elements of order 3 in the known automorphism groups: 
The group $G:=\Aut (P_{48p}) $  contains 2 conjugacy classes of elements $\sigma $ of order 3 with a 
16-dimensional fixed lattice. In both cases the automorphism group of the orthogonal lattice
has order $2^83^65$. The elements $\sigma $ may be distinguished by the order of the
normaliser $N_G(\langle \sigma \rangle )$ with is either $288$ or $144$. 
In the second case this normalizer contains a unique subgroup $U_p$ with $\sigma \in U_p$ 
such that $U_p$ is isomorphic (as an abstract group) to $U$.
We take this group $U_p$ compute the lattice $M_p$ as full preimage of the fixed space $\Fix _{U_p}(P_{48p}/2P_{48p})$.
We compute an isometry between this lattice $M_p$ and the sublattice $M$  of $A_{46}(C_{46,48})$. 
This isometry turns out to map $P_{48p}$ onto the lattice $A_{46}(C_{46,48})$.
\\
Also the group 
$G:=\Aut (P_{48n}) $  contains 2 conjugacy classes of elements $\sigma $ of order 3 with a 
16-dimensional fixed lattice. Here the automorphism groups of the orthogonal lattices
have order $2^83^9$ respectively $2^{15}3^95^2$.  The 
normaliser $N_G(\langle \sigma \rangle )$ of the first element $\sigma $ 
contains two subgroups $U_n$  with $\sigma \in U_n$ 
such that $U_n$ is isomorphic (as an abstract group) to $U$.
Only for one of these groups $U_n$ the lattice $M_n$ 
(the full preimage of the fixed space $\Fix _{U_n}(P_{48n}/2P_{48n})$) is isometric to the lattice 
$M$ constructed from $A_{14}(C_{14,48})$.
Again the isometry turns out to map $P_{48n}$ onto the lattice $A_{14}(C_{14,48})$.
\eb

\section{Automorphisms of unimodular lattices.} 

Let $\sigma \in \GL_n(\Q )$ be some element of prime order $p$. 
Let $K:=\ker (\sigma -1)$  and $I:=\im (\sigma -1)$.
Then $K$ is the fixed space of $\sigma $ and the action of 
$\sigma  $ on $I$ gives rise to a $\Q[\zeta _p]$-structure of $I$.
In particular $n= d+z(p-1)$, where 
$d:=\dim _{\Q} (K)$  and $z = \dim _{\Q[\zeta _p]} (I)$. 

If $L$ is some $\sigma $-invariant $\Z$-lattice, then 
$L$ contains a sublattice  $M$ with
$$ L \geq M:= (L\cap K) \oplus (L \cap I) = L_K(\sigma ) \oplus L_I (\sigma ) \geq pL $$ 
of finite index $[L:M] = p^s$ where $s\leq \min (d,z )$.
The fixed lattice $L_K (\sigma ) $ is sometimes also denoted by $\Fix(\sigma ) = \Fix_L(\sigma)$.

\begin{defn}
The tupel $p-(z,d)-s$ is called the {\em type} of the element $\sigma \in \GL(L)$.
\end{defn}

We now assume that we fix some $\sigma $-invariant positive definite symmetric bilinear form 
$F$. Since there are no non-zero $\sigma $-invariant homomorphisms between $K$ and $I$,
 the decomposition above is orthogonal with respect to $F$,
$\R ^n = K\perp I$. 
If we assume additionally that the $\sigma $-invariant lattice $L$ is unimodular with respect to $F$, i.e. 
$$L = L^{\#} := \{ x\in \R L \mid F(x,L) \subseteq \Z \} $$
then we obtain the following lemma.

\begin{lemma} 
Let $(L,F)$ be an even unimodular lattice and $\sigma \in \Aut (L) := \GL (L) \cap O(F) $
be of type $p-(z,d)-s$. Then  $s\leq \min (z,d )$ and
$$ L_K (\sigma )^{\#} /L_K (\sigma ) \cong L_I (\sigma )^{\#} / L_I (\sigma ) \cong (\Z/p\Z )^s $$
as abelian groups.
\end{lemma}

\bew
Put $L_K:=L_K(\sigma )$ and $L_I:=L_I(\sigma )$.
The fact that $ L_K^{\#} /L_K \cong L_I^{\#} / L_I $ is well known 
(see for instance \cite[Prop. 1.9.8]{Martinet}). \\
Let $\pi _K:= \frac{1}{p} (1+\sigma + \ldots + \sigma ^{p-1})$ denote the orthogonal projection 
onto $K$ and $\pi _I = 1 -\pi_K$ the one onto $I$. 
Then $\pi _I$ and $\pi_K$ commute with $\sigma $, $\sigma \pi _K = \pi _K$ and hence
$(1-\sigma ) \pi _I = (1-\sigma ) (1-\pi _K) = 1-\sigma $. 
The lattice $L_I = I \cap L$ contains $(1-\sigma ) L$. 
We show that $(1-\sigma )$ annihilates the quotient $L_I^{\#} / L_I$. 
Since $L=L^{\#} $ the dual lattice of $L_I$ is $L_I^{\#} = \pi _I (L)$. 
Therefore $$(1-\sigma )L_I^{\#} = (1-\sigma ) \pi _I (L) = (1-\sigma )L \subseteq L_I .$$
Therefore $L_I^{\#} / L_I $ is a quotient of  $L_I^{\#} / (1-\sigma )L_I^{\#}  \cong \F_p^z$, so 
$s\leq z$. 
Similarly $$L_K^{\#} = \pi _K (L) \supseteq L_K = \pi _K(L) \cap L \supseteq p \pi _K(L) $$
and hence $s\leq d$.
\eb

Note that for odd primes $p$ there is at most one genus of even $p$-elementary lattices of given
dimension and determinant (\cite[Theorem 13, p. 386]{SPLAG}). 
So the type of $\sigma $ 
uniquely determines  the genus of $L_K(\sigma ) $ and $L_I(\sigma ) $.

\begin{lemma} \label{kongruenz}
$z-s \in 2\Z _{\geq 0}$.
\end{lemma}

\bew
The image $I=\im (\sigma -1) $ is  
$z$-dimensional over $\Q[\zeta _p]$ and $L_I:=L_I(\sigma ) $ is a $\Z[\zeta _p]$-lattice in $I$.
Since $F$ is $\sigma $-invariant, there is some
 positive definite Hermitian  form $h:I\times I \to \Q[\zeta _p]$ so that 
$F(x,y) = \trace_{\Q[\zeta _p]/\Q } (h(x,y) ) $ for all $x,y \in I$.
Choose
some orthogonal basis $(b_1,\ldots , b_z)$, so that $h(b_i,b_j) = a_i \delta_{ij}$ 
with $a_i \in \Z[\zeta _p + \zeta _p^{-1} ]$ and put 
${\mathcal M} := \langle b_1,\ldots , b_z \rangle _{\Z [\zeta _p]} $ the $\Z[\zeta _p]$-lattice
spanned by this basis.
Let $M$ be  ${\mathcal M}$ viewed as a $\Z $-lattice together 
with the positive definite symmetric bilinear form $F$.  
Then the $\Z$-dual lattice $M^{\#} $ with respect to $F$ is
$$M^{\#} = \langle \frac{1}{a_1(1-\zeta _p)^{p-2}} b_1,\ldots , \frac{1}{a_z (1-\zeta _p)^{p-2} } b_z \rangle _{\Z [\zeta _p]}$$
and in particular $\det(M,F) = |M^{\#}/M| = p^{z(p-2)} \prod _{i=1}^z N_{\Q[\zeta _p]/\Q } (a_i) $.
Since $N_{\Q[\zeta _p]/\Q } (a_i) = N_{\Q[\zeta _p + \zeta _p^{-1}]/\Q } (a_i)^2$ 
is a rational square, 
the square class of the determinant of $F$ is $p^z (\Q ^*)^2$. 
In particular $\det (L_I) p^{z}  = p^{s+z} $ is a square so $s\equiv z \pmod{2}$. 
\eb

\subsection{Prime order automorphisms of  48-dimensional extremal lattices.} \label{prime}

In this section,
$L$ is always an extremal even unimodular lattice of dimension 48.
We list the possible types of automorphisms $\sigma \in \Aut(L)$ of prime order.

\begin{theorem}\label{autbig}
Let $L$ be an extremal even unimodular lattice of dimension $48$.
Assume that $\sigma \in \Aut(L)$ is of prime order $p\geq 11$.
Then $\sigma $ is of type 
$$ 47-(1,2)-1,\ 23-(2,4)-2 ,\ 13-(4,0)-0 ,\ 11-(4,8)-4 $$
and the fixed lattice $L_K(\sigma )$ is the unique extremal $p$-modular lattice 
of dimension $2,4,$ and $8$. Its  Gram matrix is $F_p$ with
$ F_{47}:= \left(\begin{array}{cc} 6 & 1 \\ 1 & 8 \end{array}\right) $ respectively $$ 
F_{23} = \left(\begin{array}{cccc} 6 & 3 & 2 & 2 \\ 3 & 6 & 2 & 0 
\\ 2 & 2 & 6 & 3 \\ 2 & 0 & 3 & 6 \end{array} \right) , \ 
F_{11} = \left( \begin{array}{cccccccc} 
 6&-2& 1& 3& 2& 0& 3&-1\\
-2& 6& 3& 1&-2& 1& 1& 0\\
 1& 3& 6& 3& 0& 3& 3& 2\\
 3& 1& 3& 6& 2& 3& 3& 2\\
 2&-2& 0& 2& 6& 3& 3& 3\\
 0& 1& 3& 3& 3& 6& 3& 3\\
 3& 1& 3& 3& 3& 3& 6& 2\\
-1& 0& 2& 2& 3& 3& 2& 6\\
\end{array}\right) 
$$
All these automorphisms occur for one of the three known extremal even unimodular lattices.
\end{theorem} 

\bew
Clearly $\varphi (p) = p-1  \leq 48$ implies that the prime $p$ is at most 47. 
We abbreviate $L_K:= L_K(\sigma )$, $L_I:=L_I(\sigma )$.
\begin{itemize}
\item[$p=47$:] Then $d=2$ and $z=1$. The fixed lattice $L_K$ is a 2-dimensional lattice of 
minimum $\geq 6$, determinant $47^s$ for $s=0$ or $s=1$. This immediately implies that $s=1$ and $L_K$ 
has Gram matrix $F_{47}$.
\item[$p=43$:] Then $z=1$ and $d=6$ and $L_K$ is a 6-dimensional lattice of minimum $\geq 6$ and determinant
dividing $43$. Since $6/(43^{1/6}) > 2/\sqrt{3} =\gamma (E_6)$ 
this contradicts the fact that the $E_6$ lattice
is the densest 6-dimensional lattice. 
\item[$p=41$:] Then $\dim (L_K) = 8, \min(L_K) \geq 6, \det (L_K) \mid 41 $ contradicts the fact that 
$E_8$ is the densest 8-dimensional lattice.
\item[$p=37$:] Here $\dim (L_K) = 12 , \min (L_K ) \geq 6 , \det (L_K ) \mid 37 $ contradicting the upper bounds
on the packing density from \cite{Elkies} given in Section \ref{Bounds}.
\item[$p=31$:] Here $\dim (L_K) = 18 , \min (L_K ) \geq 6 , \det (L_K ) \mid 31 $ contradicting the upper bounds
on the packing density from \cite{Elkies} given in Section \ref{Bounds}.
\item[$p=29$:] Here $\dim (L_K) = 20 , \min (L_K ) \geq 6 , \det (L_K ) \mid 29 $ contradicting the upper bounds
on the packing density from \cite{Elkies} given in Section \ref{Bounds}.
\item[$p=23$:] Then the type of $\sigma $ is either $23-(2,4)-s$ with $s\leq 2$ or 
$23-(1,26)-s$ with $s\leq 1$. 
The latter case easily gives a contradiction to the bounds in \cite{Elkies} given in Section \ref{Bounds}. 
In the other case, $\dim (L_K) =4$, $\min (L_K) \geq 6$, and $\det (L_K) \mid 23^2$ immediately yields that 
$\det(L_K) = 23^2$, so $L_K$ is in the genus of the even 23-modular lattices.
This genus contains a unique extremal lattice. 
\item[$p=19$:] Then the type of $\sigma $ is either $19-(2,12)-s$ with $s\leq 2$ and $L_K$ is a 12-dimensional
lattice of determinant dividing $19^2$ and minimum $6$ contradicting the bound in \cite{Elkies} or 
$\sigma $ has type $19-(1,30)-s$ with $s\leq 1$ and $L_I$ is an $18$-dimensional lattice of 
determinant dividing $19$ and minimum $6$ contradicting the bound in \cite{Elkies} given in Section \ref{Bounds}.
\item[$p=17$:] 
Since the kissing number of $L$ is $2^93^25^37\cdot 13$ which is not a multiple of 17,
the automorphism $\sigma $ can not act fixed point freely. This also follows from Lemma \ref{kongruenz}. So
the type of $\sigma $ is either $17-(2,16)-s$ with $s\leq 2$ 
or  $17-(1,32)-s$ with $s\leq 1$. 
In the first case $L_K$ and in the second case $L_I$ is a 16-dimensional
lattice of determinant dividing $17^2$ and minimum $\geq 6$ 
contradicting the bound in \cite{Elkies} given in Section \ref{Bounds}. 
\item[$p=13$:] 
 Then the type of $\sigma $ is either $13-(4,0)-0$ (which occurs as an element of $\Aut(P_{48n})$), 
or $13-(3,12)-s$ with $s\leq 3$ and $L_K$ is a 12-dimensional
lattice of determinant dividing $13^3$ and minimum $6$,
or $13-(2,24)-s$ with $s\leq 2$ and $L_K$ is a 24-dimensional lattice 
of determinant dividing $13^2$ and minimum $6$,
or $13-(1,36)-s$ with $s\leq 1$ and $L_I$ is a 12-dimensional lattice 
of determinant dividing $13$ and minimum $6$. 
The latter three cases
 contradict the bounds in \cite{Elkies} given in Section \ref{Bounds}.
\item[$p=11$:] 
 Then the type of $\sigma $ is either $11-(4,8)-s$ with $s\leq 4$,
or $11-(3,18)-s$ with $s\leq 3$ and $L_K$ is a 18-dimensional
lattice of determinant dividing $11^3$ and minimum $6$,
or $11-(2,28)-s$ with $s\leq 2$ and $L_K$ is a 22-dimensional lattice 
of determinant dividing $11^2$ and minimum $6$,
or $11-(1,38)-s$ with $s\leq 1$ and $L_I$ is a 10-dimensional lattice 
of determinant dividing $11$ and minimum $6$. 
The latter three cases
 contradict the bounds in \cite{Elkies} given in Section \ref{Bounds}.
In the first case, these bounds imply that $s=4$. Hence $L_K$ is in the genus of $11$-modular
$8$-dimensional lattices. By \cite{Scharlau} this genus contains a unique lattice of minimum $\geq 6$.
This lattice has Gram matrix $F_{11} $ as claimed. 
\eb
\end{itemize}

\begin{prop}\label{aut7}
Let $\sigma \in \Aut(L)$ be of order $7$. Then either the 
type of $\sigma $ is $7-(8,0)-0$ or the type is $7-(7,6)-5$ with fixed lattice
$L_K (\sigma ) \cong \sqrt{7} A_6^{\#} $.
\end{prop}

\bew
Clearly the type of $\sigma $ is $7-(8-a,6a)-s$ with $s\leq 8-a$.
For $a\geq 3$ one obtains a contradiction 
since the density of one of the lattices $L_K:=L_K(\sigma )$ or 
$L_I:=L_I(\sigma )$
exceeds the bounds  in \cite{Elkies}. 
For $a=2$ these bounds allow the possibility that $\det(L_K ) = 7^6$.
Then $L_K$ is a 12-dimensional even 7-elementary lattice in the genus of the 7-modular
lattices. The complete enumeration of this genus in \cite{Scharlau} has proved that
such a lattice does not exist.
In the case $a=1$ the lattice $L_K$ is of dimension 6, even, 7-elementary.
One concludes that $s=5$ and $L_K \cong \sqrt{7} A_6^{\# }$ because $A_6$ is the
only even 6-dimensional lattice of determinant 7. 
\eb

Only $P_{48n}$ contains an automorphism of order 7. This automorphism is of type 
$7-(8,0)-0$. Though I do not know an extremal even unimodular lattice with an
automorphism of type $7-(7,6)-5$, I cannot exclude this possibility. 
The automorphism group of the Leech lattice contains elements of type
$5-(5,4)-3$ where the fixed lattice is $\sqrt{5}A_4^{\#}$.

Using the bounds in \cite{Elkies} and Lemma \ref{kongruenz} we get the following possible 
types for automorphisms of order 5. 

\begin{rem}
Let $\sigma \in \Aut(L)$ be of order $5$. 
Then the type of $\sigma $ is  either 
$5-(12,0)-0$, $5-(10,8)-8$,  $5-(10,8)-6$, 
$5-(9,12)-9$, $5-(9,12)-7$, 
$5-(8,16)-8$,
or $5-(7,20)-7$. 
\end{rem}

Excluding some cases in this remark we obtain the following proposition:

\begin{prop}\label{aut5}
Let $\sigma \in \Aut(L)$ be of order 5. 
Then the type of $\sigma $ is  either 
$5-(12,0)-0$, $5-(10,8)-8$ with $L_K (\sigma )\cong \sqrt{5} E_8$, 
$5-(8,16)-8$ with $L_K (\sigma )\cong [2.\Alt_{10}]_{16} $ and $L_I(\sigma )$ is a hermitian unimodular lattice
over $\Z[\zeta _5]$ of dimension 8
or $5-(7,20)-7$. 
\end{prop}

\bew
For the type $5-(10,8)-8$, the fixed lattice is $\sqrt{5} U$ for some 8-dimensional even
unimodular lattice, so clearly $L_K(\sigma ) \cong \sqrt{5} E_8$.
No overlattice of index 5 of this lattice has minimum $\geq 6$ ($\Aut(E_8) $ has 8 orbits on the one-dimensional
subspaces of $E_8/5E_8$).
\\
For the type $5-(8,16)-8$, the fixed lattice is  in the genus of the even 5-modular 
16-dimensional lattices. By \cite[Theorem 8.1]{BaVe} this genus contains a unique lattice of minimum 6, denoted
by $[2.\Alt_{10}]_{16}$ in \cite{NePl}.
\\
To exclude the types 
$5-(9,12)-a$  for $a=9,7$, it is enough to note that  all 15 isometry classes of lattices in 
the genus of 
$\sqrt{5}(A_4\perp A_4 \perp A_4)^{\#} $ have minimum $\leq 4$.
\eb

Only the lattice $P_{48n}$ has an automorphism of order 5. This automorphism 
has type $5-(12,0)-0$.

\begin{prop}\label{aut3}
Let $\sigma \in \Aut(L)$ be of order $3$. 
Then the type of $\sigma $ is  either 
$3-(24,0)-0$, $3-(20,8)-8$ with $L_K(\sigma ) \cong \sqrt{3} E_8$, 
$3-(18,12)-12$ with $L_K(\sigma ) \cong \sqrt{3} D_{12}^+$, 
$3-(16,16)-16$ with $L_K(\sigma ) \cong \sqrt{3} (E_8\perp E_8) $ or $L_K(\sigma )\cong \sqrt{3} D_{16}^+$,
type $3-(15,18)-15$ (two possibilities for $L_I$, $L_K$ is unique),
type $3-(14,20)-14$ ($L_I$ unique), or
type $3-(13,22)-13$ ($L_I$ unique).
Only $L_K=0$ ($P_{48n}, P_{48p}$) and $L_K \cong \sqrt{3} (E_8\perp E_8) $ 
($P_{48n}, 2\times, P_{48p}, 2\times , P_{48q}$) occur for the three known lattices.
\end{prop}

\bew
Put $L_K:=L_K(\sigma )$, $L_I:=L_I(\sigma )$.
Most of the possible types are excluded from the bounds in \cite{Elkies} using Lemma \ref{kongruenz}. 
So I only comment on the cases, where these bounds do allow lattices of minimum 6. These are \\
type $3-(19,10)-9$:   Then $L_K$ is in the genus of $\sqrt{3} E_8\perp A_2 $.
This is a one-class genus, so $\min (L_K) = 2 < 6$ is a contradiction.
\\
type $3-(18,12)-10$: Then the lattice $L_K$ is an overlattice of index 3 of $\sqrt{3} D_{12}^+$.
The automorphism group of $D_{12}^+$ has 13 orbits on these lattices, non of them has minimum $6$.
\\
type $3-(17,14)-13$: Then the lattice $L_K$ is in the genus of $\sqrt{3} (E_8 \perp E_6^{\#} )$.
This genus contains 2 isometry classes, both lattices have minimum 4. 
This also excludes the type $3-(17,14)-11$, since then $L_K$ contains one of the latter two lattices. 
\\
To exclude the cases $3-(16,16)-s$ for $s=14,12,10$ we need to compute overlattices of index 3 of 
the rescaled unimodular lattices $\sqrt{3} (E_8\perp E_8) $ and $\sqrt{3} D_{16}^+$.
The automorphism group of the first lattice has 14 orbits, the other one 17 orbits on the overlattices of index 3,
none of the lattices has minimum 6. 
\\
types $3-(15,18)-a$, $a=15,13,11$;
 $3-(14,20)-a$, $a=14,12,10$;
 $3-(13,22)-a$, $a=13,11$;
 $3-(12,24)-a$, $a=12,10$;
$3-(11,26)-11$;
 $3-(10,28)-10$: 
\\
 In all cases the lattice $L_I$ is or contains  an Hermitian unimodular lattice
of dimension $z=10,\ldots , 15$. These lattices have been classified by Feit, Abdukhalikov, and Scharlau
\cite{F}, \cite{A}, \cite{AS}. 
For lattices of minimum 6 one obtains
one lattice, $L_{13}$, of dimension $2z=26$, 
one lattice, $L_{14}$, of dimension 
$2z=28$, 
and two such lattices, $L_{15}$, $M_{15}$, of dimension $2z=30$.
\\
This excludes the cases $z=10,11,12$. 
\\
The group $\Aut(L_{13})$ has 2 orbits on the set of integral 
overlattices  of index 3 of $L_{13}$, one has minimum 2 and the other minimum 4. 
So $L_I=L_{13}$ and only the type $3-(13,22)-13$ is possible here.
\\
The group $\Aut(L_{14})$ has 5 orbits on the set of integral 
overlattices of index 3 of $L_{14}$, one with minimum 2, the other 4 with minimum 4.
So $L_I=L_{14}$ and only the type $3-(14,22)-14$ is possible here.
\\
The group $\Aut(L_{15}) \cong \pm 3.U_4(3).2$ has 13 orbits on the set of integral
 overlattices of index 3 of $L_{15}$,
again none of these overlattices has minimum $\geq 6$.
\\
The group $\Aut(M_{15}) \cong \pm (3^{1+2}_+ \times 3^{1+2}_+) .\SL_2(3).2$ has 
174 orbits on the set of integral overlattices  of index 3 of $M_{15}$,
again none of these overlattices has minimum $\geq 6$.
\eb

\begin{lemma}\label{2elt}
Let $M$ be an even lattice such that $M^{\#}/M$ has exponent $2$.
Then $M$ contains a sublattice isometric to $\sqrt{2} U$ where 
$U=U^{\#}$ is an (even or odd) unimodular lattice.
\end{lemma}

\bew
Since this is a statement about 2-adic lattices, we pass to $M_2:=\Z_2\otimes M$. 
This lattice has a 2-adic Jordan decomposition $M_2 = f_{II} \oplus \sqrt{2} f $ where
$f_{II}$ is even and unimodular of dimension, say, $2m$, and $f$ is unimodular. 
If $m=0$, then $M \cong \sqrt{2} U$ for some unimodular $U$ and we are done.
So assume $m\geq 1$. Then
$f_{II}$ contains a vector $v$ such that $(v,v) \in 2 \Z_2^{*}$. 
Therefore $f_{II}$ contains the sublattice $v \perp v^{\perp} $ of index 2, and 
$v \perp v^{\perp} \perp \sqrt{2} f $ has a Jordan decomposition $g_{II} \oplus \sqrt{2} g$ 
with $\dim (g) = \dim (f) +2$ and $\dim (g_{II} ) = 2(m-1)$. 
Since $g_{II}$ is again even and unimodular, we may proceed by induction.
\\
Note that 
the possible genera of such lattices are given in \cite[Table 15.5, p. 388]{SPLAG} and
one can easily find lattices representing these genera and give a case by case proof of the lemma.
\eb

\begin{theorem}\label{aut2}
Let $-1\neq \sigma \in \Aut(L) $ be of order $2$. 
Then $\sigma $ is of  type $2-(24,24)-24$ and 
$$\Fix _L(\sigma ) \perp \Fix_L(-\sigma )  
\cong \sqrt{2} (\Lambda _{24} \perp \Lambda _{24}) \mbox{ or } 
 \sqrt{2} (O_{24} \perp O_{24}) .$$
  Both cases occur.
\end{theorem}

\bew 
Both lattices $L_K:=L_K (\sigma ) $ and $L_I:=L_K(-\sigma )$ satisfy that 
$$ L_K^{\#}/L_K \cong L_I^{\#}/L_I \mbox{ has exponent } 2, \min (L_K ) \geq 6, \min (L_I) \geq 6, 
L_I \mbox{ and } L_K \mbox{ even.} $$
By Lemma \ref{2elt} 
all such lattices contain some lattice $\sqrt{2} U$ with $U=U^{\#}$. 
Since there is no unimodular lattice $U$ of dimension $n\leq 22$ or of dimension 25  
(see \cite[Table 16.7, p. 416-417]{SPLAG} and \cite{Borcherds})
with $\min (U) \geq 3$ this implies that $\dim (L_K ) = 24$  and that 
$L_K $ is an overlattice of either $\sqrt{2}\Lambda _{24} $ or $\sqrt{2} O_{24}$
where $\Lambda _{24}$ and $O_{24}$ denote the Leech lattice respectively the odd Leech lattice
of dimension 24 (see \cite[Table 16.1, p. 407]{SPLAG} and \cite[Table 17.1, p. 424-426]{SPLAG}). 
All non-zero classes of $\frac{1}{\sqrt{2}}\Lambda _{24}/\sqrt{2}\Lambda _{24}$ are 
represented by vectors of norm 2,3,4, (see e.g. \cite[xxx]{Ebeling}) so $\sqrt{2}\Lambda _{24}$  has no overlattice of index
2 with minimum $\geq 6$. 
The automorphism group of the odd Leech lattice has 16 orbits on the 1-dimensional subspaces
of $\frac{1}{\sqrt{2}} O_{24}/\sqrt{2} O_{24}$. The minima of the corresponding overlattices are
$$\frac{3}{2}, 2 \mbox{ (2 lattices)}, 
\frac{5}{2}, 3 \mbox{ (2 lattices)}, 
\frac{7}{2} \mbox{ (2 lattices)}, 4 \mbox{ (5 lattices)}, 
\frac{9}{2} \mbox{ (2 lattices)}, 5  .$$
In particular no proper overlattice of index 2 has minimum $\geq 6$.
\\
The case $\Fix _L(\sigma ) \perp \Fix_L(-\sigma )  
\cong \sqrt{2} (O _{24} \perp \Lambda _{24}) $ cannot occur, since the quadratic spaces
$\frac{1}{\sqrt{2}} O_{24}/\sqrt{2} O_{24} $ and 
$\frac{1}{\sqrt{2}} \Lambda _{24}/\sqrt{2} \Lambda _{24} $ are not isometric.
\eb

\begin{center}
The following table summarizes the results of this section. 
\\
\begin{tabular}{|c|c|c|c|c|}
\hline
p & $\dim L_K (\sigma ) $ & $\det (L_K (\sigma )) $ & $L_K  (\sigma ) $ & example \\
\hline
47 & 2 & 47 & unique & $P_{48q}$ \\ 
\hline
23 & 4 & $23^2$  & unique & $P_{48q}$, $P_{48p}$ \\
\hline
13 & 0 &   & $\{ 0 \}$ & $P_{48n}$ \\
\hline
11 & 8 & $11^4$  & unique & $P_{48p}$ \\
\hline
7 & 0 &   & $\{ 0 \}$ & $P_{48n}$ \\
7 & 6 & $7^5$  & $\sqrt{7}A_6^{\#}$ & no \\
\hline
5 & 0 &   & $\{0 \} $ & $P_{48n}$ \\
5 & 8 & $5^8$  & $\sqrt{5}E_8$ & no \\
5 & 16 & $5^8$  & $[2.\Alt_{10}]_{16}$ & no \\
5 & 20 & $5^7$ & ? & no \\
\hline
3 & 0 & & $\{0 \}$  & $P_{48p}$, $P_{48n}$ \\
3 & 8 & $3^8$ & $\sqrt{3} E_8$ & no \\
3 & 12 & $3^{12}$ & $\sqrt{3} D_{12}^+$ & no \\
3 & 16 & $3^{16}$ & $\sqrt{3} (E_8\perp E_8) $ &  $P_{48p}$, $P_{48q}$, $P_{48n}$  \\
3 & 16 & $3^{16}$ & $\sqrt{3} D_{16}^+$ & no \\
3 & 18 & $3^{15}$ & unique & no \\
3 & 20 & $3^{14}$ & $L_I (\sigma )$ unique & no \\
3 & 22 & $3^{13}$ & $L_I (\sigma )$ unique & no \\
\hline
2 & 0 &   & $\{ 0 \}$ & $\sigma = -1$ \\
2 & 24 & $2^{24}$  & $\sqrt{2}\Lambda _{24} $ & $P_{48n}$ \\
2 & 24 & $2^{24}$  & $\sqrt{2} O _{24} $ & $P_{48n}$, $P_{48p}$ \\
\hline
\end{tabular} 
\end{center}

\begin{rem}
$\Aut (P_{48n})$ has $3$ conjugacy classes of elements $\sigma \neq -1$ of order $2$,
two of which have $L_K(\sigma ) \cong L_K(-\sigma ) \cong \sqrt{2}\Lambda _{24} $ 
and one has $L_K(\sigma ) \cong L_K(-\sigma ) \cong \sqrt{2}O _{24} $.
Both classes of elements $\sigma \neq -1$ of order $2$ in $\Aut(P_{48p})$ have 
 $L_K(\sigma ) \cong L_K(-\sigma ) \cong \sqrt{2}O _{24} $ whereas $\Aut (P_{48q})$ 
does not contains nontrivial elements of order $2$.
\end{rem}

\begin{theorem}
Let $L$ be an extremal even unimodular lattice of dimension $48$ such that $\Aut(L) $ contains 
an element $\sigma  $ of order $46$ such that $\sigma ^{23} \neq -1$.
Then $L\cong P_{48p}$.
\end{theorem}  

\bew
The automorphism group of the Leech lattice $\Aut(\Lambda _{24}) = 2.Co_1$ and 
the odd Leech lattice $\Aut(O_{24}) \cong 2^{12} . M_{24} $ both contain two 
conjugacy classes of elements of order 23 represented by $g$ and $g^{-1}$, say.
A  Gram matrix for the respective fixed lattices is given by 
$$ F_{\Lambda} := \left( \begin{array}{cc} 4 & 1 \\ 1 & 6 \end{array}\right) 
\mbox{ and } 
F_O := \left( \begin{array}{cc} 3 & 1 \\ 1 & 8 \end{array}\right)  .$$
The automorphism group of the 4-dimensional lattice $X$ with Gram matrix $F_{23}$ 
of Theorem \ref{autbig} has 2 conjugacy classes of automorphisms $\neq -1$ of order 2,
represented by, say, $h_1$ and $h_2$. 
Since $h_i$ and $-h_i$ are conjugate in $\Aut(X)$ the fixed lattice of $h_i$ and $-h_i$ are isometric,
for $h_1$ they have Gram matrix $2F_O$ and for $h_2$ the Gram matrix is $2F_{\Lambda }$.
\\
Let $L_K:=L_K(\sigma ^{23})$ and $L_I:=L_K(-\sigma ^{23})$.
From  Theorem \ref{aut2} one concludes that $L_K \cong L_I \cong \sqrt{2} \Lambda _{24}$ or
 $L_K \cong L_I \cong \sqrt{2} O _{24}$. 
By Theorem \ref{autbig}
the element $\sigma ^{2}$ acts on $L_K\perp L_I$ as an automorphism of
 type  $23-(2,4)-2$.
The conjugacy classes of such subgroups $\langle \sigma ^2 \rangle $ of order 23 in 
$\Aut (L_K \perp L_I)$ are represented by $\langle (g,g)\rangle $ and $\langle (g,g^{-1}) \rangle$.
With MAGMA we compute the $\langle \sigma ^2 \rangle $ invariant unimodular overlattices 
$L$ of $L_K\perp L_I$
of minimum 6 for all 4 cases. 
Note that $\frac{1}{2} (L_K\perp L_I)/(L_K\perp L_I) $ is a semisimple $\F_2 \langle \sigma ^2 \rangle $-module
isomorphic to $S^4\oplus V_1^2 \oplus V_2^2$ where $S \cong \F_2$ and $V_1$, $V_2 \cong \Hom _{\F_2}(V_1,\F_2)$
 are the two non-isomorphic 11-dimensional $\F_2 \langle \sigma ^2 \rangle $-modules.
Then $L/(L_K\perp L_I) \cong S^2 \oplus V_1 \oplus V_2$ can be reached in 2 steps: 
First we compute the unique lattice $Y$ of minimum 6 with $Y/(L_K\perp L_I) \cong S^2$,
by going through all 35 two-dimensional subspaces of $S^4 \cong \F_2^4$. 
We then go through the $(2^{11}-1)/23 +2 = 91$  orbits of one-dimensional subspaces 
of $V_1 \perp V_1$ under the action of $\langle (1,g)\rangle \leq C_{\Aut(L_K \perp L_I)} (\sigma ^2)$
to obtain candidates for the lattices $Z$ of minimum 6 such that 
$Z/Y \cong V_1$. Then $Z^{\#}/Z \cong V_1 \oplus V_2 $ and there is a unique unimodular 
overlattice $W$ of $Z$ with $W/Z \cong V_2$. 
Only for the case $L_K \cong \sqrt{2} O_{24}$  and 
$\langle \sigma ^2 \rangle \cong \langle (g,g^{-1}) \rangle$ there is such a lattice $W$ 
with $\min (W) =6 $. We check that $W\cong P_{48p}$ by 
computing a vector $\beta \in W$ fixed by $\sigma ^2$ and of norm $(\beta,\beta ) = 12$ 
such that the neighbor $W^{(\beta )}$ contains a 3-frame. 
The corresponding extremal ternary code is easily checked to be isometric to the 
Pless code. 
\eb

From the previous discussion we found that the maximal dimension
of the fixed lattice of an  automorphism of $L$ of odd prime order  is $\leq 22$.
The nontrivial automorphism of order 2 have fixed lattices of dimension $24$. 
From this fact we obtain the following corollary.

\begin{cor}\label{Phi}
Let $L$ be an extremal even unimodular lattice of dimension $48$
and $\sigma \in \Aut(L) $ of order $m$.
Then the $m$-th cyclotomic polynomial $\Phi _m$ divides the
minimal polynomial of $\sigma $. 
\end{cor}

\bew
Assume that $\Phi _m $ does not divide $ \mu _{\sigma }$. Then there are two distinct 
prime divisors $p$ and $q$ of $m$ such that any element $a
\in \langle \sigma \rangle $  of order $pq$ has minimal polynomial dividing $\Phi _p \Phi _q \Phi _1$. 
This implies that $\Fix (a^p) + \Fix  (a^q)$ generates $\R ^{48}$. 
By the results in this section
 the maximal dimension of a fixed space of some non-trivial element in 
$\Aut (L)$ is 24, and this only occurs for elements of order 2, 
for bigger primes, the dimension of the fixed space is at most 22, so the sum
 $\Fix (a^p) + \Fix  (a^q)$ has at most dimension $46$, a contradiction. 
\eb

\begin{cor}\label{cor2}
Let $L$ be an extremal even unimodular lattice of dimension $48$.
Then any element $\sigma $ of order $2^5$ in 
$\Aut (L) $ acts with irreducible minimal polynomial 
and $\Aut(L)$ contains no elements of order $2^6$, $2^53$,
 $2^45$, or 
$2^3 11$.
\end{cor}

\bew
First assume that $\sigma \in \Aut(L)$ has order $2^5$. 
Since $\varphi (2^5) = 16$ and also $\dim (L) = 48$ is a multiple of
16, also the dimension of the fixed lattice 
$\Fix _L(\sigma ^{2^4})$ is a multiple of 16. 
By Theorem \ref{aut2} this is only possible, if $\Fix _L(\sigma ^{2^4}) = \{ 0 \} $ so $\sigma ^{2^4} = -1$.
\\
Now assume that $\sigma \in \Aut(L)$ has order $2^6$. 
Then $\sigma ^2$ acts with irreducible minimal polynomial,
contradicting the fact that 
 $\varphi (2^6) = 32$ does not divide $\dim (L) = 48$.
\\
Now let $|\langle \sigma \rangle | = 2^5 3$. 
Then $\sigma ^3$ acts with an irreducible minimal polynomial on
$L$, so the type of $\sigma ^{2^5}$ is $3-(z,d)-s$ with 
$z$ and $d$ both divisible by 16. 
By Proposition \ref{aut3} this yields that 
this type is $3-(16,16)-16$ and 
$\Fix_L(\sigma ^{2^5}) $ is either $\sqrt{3} (E_8\perp E_8)$ or
$ \sqrt{3} D_{16}^+$. With Magma we check that both lattices admit
no automorphism of order 32.
\\
Now let $|\langle \sigma \rangle | = 2^4 5$ and put $\tau := \sigma ^{2^35} $.
Then $\tau $ is an element of order 2 in $\Aut (L)$.
If $\tau\neq -1$ then by Theorem \ref{aut2}
both lattices $\Fix_L(\tau )$ and 
$\Fix_L(-\tau )$ are either similar to $\Lambda _{24}$ or
to $O_{24}$. 
By Proposition \ref{aut5} one sees that $\sigma $  
acts as an element of order 80 on $\Fix _L(-\tau )$.
Neither $\Aut(\Lambda _{24})$ nor $\Aut (O_{24})$ contain 
elements of order 80.
\\
Therefore  $\tau = -1$, so $\sigma ^5$ acts with
irreducible minimal polynomial.  As above this implies that 
the type of $\sigma ^{2^4}$ is $5-(z,d)-s$ with $z$ and $d$ 
both divisible by $8$. By Proposition \ref{aut5} this implies that 
$\Fix _L(\sigma ^{2^4}) =L_K \cong [2.\Alt_{10}]_{16}$. 
With Magma one checks that $\Aut(L_K) \cong 2.\Alt_{10} $ does not contain an 
element of order 16. 
\\
Finally let $\sigma $ be of order $2^311$, $L_K:=L_K(\sigma ^8)$ 
and $L_I:=L_I(\sigma ^8)$.
Since $\Aut (L_K)$ does not contain an element of order 8,
the element $\sigma $ acts with an irreducible minimal polynomial 
on $L_I$ (which is hence an ideal lattice in the 88-th cyclotomic
number field). But then $L_K = L_K (\sigma ^{44})$
 contradicting Theorem \ref{aut2} 
(and also the fact that $L_K$ has determinant $11^4$). 
\eb

\section{Lattices with a given automorphism}

In this section we compute all extremal even unimodular lattices $L$ of dimension 48 that
admit an automorphism $\sigma $ of order $m$ with $\varphi (m) > 24$. 
This can be reduced to 
a computation with ideals in the $m$-th cyclotomic number field.
We first deal with the easier case $\varphi (m) = 48$, where $L$ is 
an ideal lattice in $\Q[\zeta _m ]$. 

\subsection{Ideal lattices} 

This section classifies all extremal even unimodular lattice $L$ that have an 
automorphism $\sigma $ of order $m$ such that $\varphi(m) = 48$. 
By Corollary \ref{Phi} the minimal polynomial of $\sigma $ is the $m$-th cyclotomic 
polynomial and hence $L$ is an ideal lattice in the $m$-th cyclotomic number field $F = \Q[\zeta _m]$.
By \cite{Bayer} for all these fields there is some positive definite even unimodular ideal lattice. 

\begin{rem} (see \cite{ideal}) 
Let $\sigma  \in \Aut (L)$ be an automorphism of the lattice $L$ with characteristic polynomial 
$\Phi _m$, the $m$-th cyclotomic polynomial.
Then the action of $\sigma $ on $\Q L$ turns the vector space $\Q L$ into a one-dimensional 
vector space over 
the $m$-th cyclotomic number field $F = \Q[\zeta _m]$ which we identify with $F$.
Then the lattice $L$ is a $\Z [\zeta _m]$-submodule, hence isomorphic to a
fractional  ideal  $J$ in $F$. 
The symmetric positive definite bilinear form $B:L\times L \to \Q $ is $\zeta _m$-invariant,
since $B(x\sigma , y \sigma ) = B(x,y) $ for all $x,y \in L$. It hence corresponds to some
trace form on the ideal $J$, 
$$B(x,y) = \trace _{F/\Q} (\alpha x \overline{y} ) = b_{\alpha }(x,y) ,\ \mbox{ for all } x,y\in J  $$
where $\overline{\phantom{x}}$ is the complex conjugation on $F$, the involution 
with fixed field $F^+:= \Q [\zeta _m + \zeta _m^+]$, and $\alpha \in F^+$ 
is totally positive, i.e. $\iota (\alpha ) > 0$ for all embeddings $\iota : F^+ \to \R $. 
Let $$\Delta := \{ x\in F \mid \trace _{F/\Q } (x \overline{y} ) \in \Z \mbox{ for all } 
y\in \Z [\zeta _m ] \} $$
denote the different of the ring of integers $\Z _F = \Z [\zeta _m]$ of  $F$. 
Then the dual lattice of $(J,b_{\alpha }) $ is 
$$(J,b_{\alpha })^{\#} = (\overline{J}^{-1} \Delta \alpha ^{-1} ,b_{\alpha }) .$$
So $(J,b_{\alpha })$ is unimodular, if and only if $(J\overline{J})^{-1} \Delta \alpha ^{-1} = \Z_F$.
For all fields $F=\Q[\zeta _m]$ it turns out that that there is some 
totally positive $\alpha_0 \in F^+$ such that $\Delta \alpha_0 ^{-1} = \Z _F$, 
so $(\Z_F,b_{\alpha _0})$ is unimodular. 
Then the isometry classes of unimodular (positive definite) ideal lattices are
represented by 
$$( J, b_{\alpha _J u \alpha _0} ) $$
where $ [J]  $ runs through all ideal classes in $F$ such that 
$J\overline{J} \cap F^+ = (\alpha _J)$ is a principal ideal generated by some totally positive element
$\alpha _J \in F^+$ and 
$$ [u] \in \{ u \in \Z_{F^+}^* \mid u \mbox{ totally positive } \} / \{ v\overline{v} \mid 
v \in \Z_F^* \} .$$
\end{rem}

\begin{definition}
Let  $F= \Q[\zeta _m]$ with maximal real subfield $F^{+} := \Q[\zeta _m + \zeta _m^{-1}]$. 
Let $\Cl(F)$ denote the ideal class group of $F$ and $\Cl^+(F^+)$ the ray class group of $F^+$,
where two ideals $I,J $ are equal in $\Cl^+(F^+)$ if there is some totally positive $\alpha \in F^+$
with $I=\alpha J$. Since $x\overline{x} $ is totally positive for any $0\neq x\in F$, 
the norm induces a group homomorphism 
$$N: \Cl(F) \to \Cl^+(F^+), \ [J] \mapsto [J\overline{J} \cap F^+]. $$
The {\bf positive class group}  $\Cl^+(F)$ is the kernel of this homomorphism and $h^+(F) := | \Cl^+(F) | $. 
\end{definition}

In the appendix of 
\cite{Washington}  the candidates for $m$ are listed as well as the class numbers of $F$,
which are correct under the generalized Riemann hypothesis.
The following table gives the structure of the class group and 
the positive class group as computed with Magma (also assuming GRH).

$$
\begin{array}{|l|c|c|c|c|c|} 
\hline
m & 65 & 104 & 105 & 112 & 140  \\
h_F & 2\cdot 2 \cdot 4 \cdot 4 & 3\cdot 117  & 13 & 3\cdot 156 & 39 \\
h^+_F & 2\cdot 2 & 3\cdot 117 &  13 & 6\cdot 39 & 39 \\
\hline
m & 144 & 156 & 168  & 180  & \\
h_F  & 13 \cdot 39 & 156 & 84  & 5\cdot 15 & \\
h^+_F  & 13 \cdot 39 & 78 & 42  & 5\cdot 15 & \\
\hline
\end{array}
$$

\begin{theorem}\label{ideal} 
Let $L$ be an extremal even unimodular lattice such that $\Aut (L)$ contains some 
element $\sigma $ of order $m$ with $\varphi (m) = 48$. 
Then $m=65$  or $m=104$ and $L\cong P_{48n}$.
\end{theorem}

\bew
We briefly describe the Magma computations that led to this result.
They are similar for all cases.
For all fields $F= \Q[\zeta _m]$ the maximal real subfield $F^{+} := \Q[\zeta _m + \zeta _m^{-1}]$
has class number 1. 
In particular the positive class group $\Cl^+(F) $ contains $\Cl(F)^2$.
In all cases there is  some totally positive $\alpha \in F^+$ such that 
$(\Z_F, b_{\alpha } )$ is unimodular. 
Let $U$ denote a system of representatives of totally positive units $u\in \Z_{F^+}^*$ modulo 
squares. 
Let $I_1,\ldots , I_k$ be ideals of $\Z_F$ such that their classes generate the positive class group of $F$,
so that $[I_j]$ has order $a_j$ with $a_1\cdots a_k = h^+_F$. 
Then there are totally positive 
$\alpha _1,\ldots , \alpha _k \in F^+$ such that $N(I_j) = (\alpha _j) $. 
Then the isometry classes of unimodular ideal lattices are all represented by some ideal lattice of the form
$$( I_1^{i_1} \cdots I_k^{i_k} , b_{u \alpha \alpha _1^{-i_1} \cdots  \alpha _k^{-i_k}} )  \mbox{ where } 
u \in U \mbox{ and } 0\leq i_l \leq a_l-1 \mbox{ for } l=1,\ldots k  .$$

\noindent
{\bf m=65:} Here $|U| = 32$, $\Cl^+(F) = \Cl(F)^2  = \langle [I_1],[I_2] \rangle  \cong C_2\times C_2 $,
$\alpha = \frac{1}{65} \alpha _{5a} \alpha _{5b} \alpha _{5c} \alpha _{13}$, so there are $4\cdot 32$ ideal 
lattices to be considered. 
Only two of them are extremal, both are principal ideal lattices. 
To obtain an isometry with the lattice $P_{48n}$ we compute a lattice $M$ such that 
$M/13L$ is the (4-dimensional) fixed space of the action of $\sigma ^{5}$ on $L/13L$. 
The dual $D$ of $M$ (rescaled to be integral) is a lattice of determinant $13^8$, minimum 6 and 
kissing number 6240 and automorphism group of order $2^6  3^2 65$. 
The automorphism group of $P_{48n}$ contains a unique conjugacy class of elements of order 65. 
We compute the corresponding sublattice $D_0$  of $P_{48n}$ and an isometry between $D_0$ and $D$. 
For both ideal lattices the isometry maps the ideal lattice to $P_{48n}$.
\\
{\bf m=104:}  Now $|U|=2$, $\Cl^+(F) = \Cl(F)  = \langle [I_1] , [I_2] \rangle  \cong C_{3} \times C_{117} $, 
$\alpha = \frac{1}{52} \alpha  _{13}$, so there are $2\cdot 3 \cdot 117 $ ideal
lattices to be considered.
Only four of them are extremal, the underlying ideal is a suitable ideal of order 3 and its inverse
and all $u\in U$.
To obtain an isometry with the lattice $P_{48n}$ we compute a lattice $M$ such that 
$M/13L$ is the (4-dimensional) fixed space of the action of $\sigma ^{8}$ on $L/13L$. As in the case $m=65$ 
the dual $D$ of $M$ (rescaled to be integral)
is isometric to the corresponding sublattice $D_0$  of $P_{48n}$ and an isometry between $D_0$ and $D$ 
maps the ideal lattice $L$ to $P_{48n}$ in all four cases.
\\
{\bf m=105:} Again $|U|=2$, $\Cl^+(F) = \Cl(F)  = \langle [I_1] \rangle  \cong C_{13} $, 
where $I_1$ is some prime ideal dividing 29.
$\alpha = \frac{1}{105} \alpha _3 \alpha _{5a} \alpha _{5b} \alpha _7$, so there are $2\cdot 13$ ideal
lattices to be considered.
In all these lattices random reduction algorithms find elements of norm 4. 
\\
{\bf m=112:} Then $|U|=4$, $\Cl^+(F) = \Cl(F)^2  = \langle [I_1],[I_2] \rangle  \cong C_6\times C_{39} $, 
where $I_1$ is some prime ideal dividing 7 and $I_2$ an ideal dividing 113. 
$\alpha = \frac{1}{56} \alpha _{7a} \alpha _{7b}$, so there are $4\cdot 6\cdot 39$ ideal
lattices to be considered.
In all these lattices random reduction algorithms find elements of norm 4. 
\\
{\bf m=140:} Here $|U|=2$, $\Cl^+(F) = \Cl(F)  = \langle [I_1] \rangle  \cong C_{39} $. 
$\alpha = \frac{1}{70}  \alpha _{5} \alpha _7$, so there are $2\cdot 39$ ideal
lattices to be considered.
In all these lattices random reduction algorithms find elements of norm 4. 
\\
{\bf m=144:} Again $|U|=2$, $\Cl^+(F) = \Cl(F)^2  = \langle [I_1],[I_2] \rangle  \cong C_{13}\times C_{39} $, 
$\alpha = \frac{1}{72} \alpha _{3} ^3$, so there are $2\cdot 13\cdot 39$ ideal
lattices to be considered.
In all these lattices random reduction algorithms find elements of norm 4. 
\\
{\bf m=156:} Now $|U|=4$, $\Cl^+(F) = \Cl(F)^2  = \langle [I_1]^2 \rangle  \cong C_{78} $, 
$\alpha = \frac{1}{78} \alpha _{3a} \alpha _{3b} \alpha _{13a} \alpha _{13b}$, so there are $4\cdot 78$ ideal
lattices to be considered.
In all these lattices random reduction algorithms find elements of norm 4. 
\\
{\bf m=168:} Again $|U|=4$, $\Cl^+(F) = \Cl(F)^2  = \langle [I_1]^2 \rangle  \cong C_{42} $, 
$\alpha = \frac{1}{84} \alpha _{3a} \alpha _{3b} \alpha _{7a} \alpha _{7b}$, so there are $4\cdot 42$ ideal
lattices to be considered.
In all these lattices random reduction algorithms find elements of norm 4. 
\\
{\bf m=180:} Now $|U|=2$, $\Cl^+(F) = \Cl(F)  = \langle [I_1] , [I_2] \rangle  \cong C_5\times C_{15} $. 
$\alpha = \frac{1}{90}  \alpha _{3}^3 \alpha _5$, so there are $2\cdot 5 \cdot 15$ ideal
lattices to be considered.
In all these lattices random reduction algorithms find elements of norm 4. 
\eb

I thank Claus Fieker for helping me with the Magma calculations.

\subsection{Subideal lattices} 

This section classifies all extremal even unimodular lattice $L$ that have an 
automorphism $\sigma $ of order $m$ such that $24 < \varphi(m) < 48$. 
By Corollary \ref{Phi} we know that $\Phi _m$ divides $\mu _{\sigma }$.
We have a unique decomposition
$$\Q L = V \oplus W \mbox{ into } \sigma \mbox{ invariant subspaces } $$
so that 
$\sigma $ acts as a primitive $m$-th
root of unity on $V$. Since $\varphi(m) > 24$, the action of $\sigma $ turns $V$ into a 1-dimensional 
vector space over $F:= \Q[\zeta _m ]$, the $m$-th cyclotomic number field, 
and the lattice $M:=L\cap V = (J,b_{\alpha })$ is an ideal lattice in $F$. 
  The lattice $M$ is integral, if $J \overline{J} \alpha \subseteq \Delta $ and then 
$$M^{\#}/M \cong \Delta (J\overline{J} \alpha ) ^{-1} /\Z[\zeta _m] \mbox{ as a } \Z[\zeta _m]-\mbox{module} .$$
For our computations it turns out that $M= L_I(\sigma ^d)$ 
for some element $\sigma ^{d}$ of prime order $p$ with $pd=m$, so we know $\det (M)$
(more precisely the abelian group $M^{\#}/M$) from the previous section. 
We also know the fixed lattice $K=\Fix_L (\sigma ^d)$ and the possible actions of $\sigma $ 
on $K$, by computing the conjugacy classes of automorphisms of $K$ of order $d$. 
The even unimodular lattice  $L$ is a subdirect product of $M^{\#} $ and $K^{\#}$ 
with kernel $M\perp K$. Therefore the $\Z[\zeta _m]$-module 
$M^{\#}/M$ is isomorphic to $K^{\#}/K$. 

\begin{center}
{\bf The values $m$ with $m$ odd or divisible by 4, such that $48 > \varphi (m) > 24 $:}  \\
\begin{tabular}{|c|c|c|c|c|c|c|c|c|c|c|c|c|}
\hline 
m & 29 & 31 & 37 & 41 & 43 & 47 & 49 & 51 & 55 & 57 & 63 & 64 \\
$\varphi(m) $ & 28 & 30 & 36 & 40 & 42 & 46 & 42 & 32 & 40 & 36 & 46 & 32 \\
thm & \ref{autbig} & \ref{autbig} &\ref{autbig} &\ref{autbig} &\ref{autbig} &\ref{aut47} & \ref{varphikleiner48} &
\ref{autbig} & \ref{aut55} & \ref{autbig} & \ref{varphikleiner48} & \ref{cor2} \\
\hline
 m & 68 & 69 & 75 & 76 & 80 & 88 & 92 & 96 & 100 & 108 & 120 & 132 \\
$\varphi(m)$ & 32 & 44 & 40 & 36 & 32 & 40 & 44 & 32 & 40 & 36 & 32 & 40 \\
thm & \ref{autbig} & \ref{varphikleiner48} & \ref{varphikleiner48} & \ref{autbig} & \ref{cor2} & \ref{cor2} & \ref{varphikleiner48}  &
\ref{cor2} & \ref{varphikleiner48} & \ref{varphikleiner48} & \ref{varphikleiner48} & \ref{varphikleiner48}  \\
\hline
\end{tabular}
\end{center}

To explain the strategy we will give two proofs in detail. 

\begin{theorem}\label{aut47}
Let $L$ be an extremal even unimodular lattice of dimension $48$ such that $\Aut(L) $ contains an
element $\sigma $ of order $47$.
Then $L\cong P_{48q}$.
\end{theorem} 

\bew
Let $L_K:= L_K(\sigma )$ and $L_I :=L_I(\sigma ) $.
Then $L_K$ has Gram matrix $F_{47}$ from Theorem \ref{autbig} and $[L:L_K\perp L_I] = 47$. 
Moreover $L_I$ is an ideal lattice in the 47-th cyclotomic number field $F=\Q[\zeta _{47}]$.
So there is some fractional ideal $J $ in $F$ and some totally positive $a\in F^+:=\Q[\zeta _{47} + \zeta _{47}^{-1}]$
such that $L_I = J$ and $(x,y) = \trace (ax\overline{y}) $ for all $x,y\in J$.
The class number of $F$ is $695=5\cdot 139$ (\cite{Washington}) and one computes with Magma, that the class group is 
generated by any prime ideal $\wp $ that divides 283 (e.g. $\wp = (283, \zeta _{47} +279 ) $). 
Again with Magma one computes a system of independent units in the ring of integers 
$\Z _{F^+} = \Z[\zeta _{47} + \zeta _{47}^{-1}]$ of which the 23 real embeddings (together with $-1$) 
generate all $2^{23} $ possible sign combinations. In particular any totally positive unit 
in $\Z _{F^+} $ is a square. Since $F^+$ has class number one, any ideal in $\Z_{F^+}$ has 
some totally positive generator. 
In particular the ideal $\wp \overline{\wp} \cap \Z _{F^+}$ has some totally  positive generator  $\alpha $
that may be computed explicitly in Magma.
Therefore the lattice $L_I$ is isometric to one of the 695 ideal lattices 
$$J_j:=(\Delta \wp ^{j} , (x,y) \mapsto 47 \trace (\alpha ^{-j} x \overline{y}))  \mbox{ for some } j = 0,\ldots , 694 .$$
Here $\Delta = (1-\zeta _{47})^{-45} $ is the different of the lattice $\Z _F$, 
which has the desired property that $\det (\Delta , (x,y)\mapsto 47 \trace(x\overline{y})) =47 $.
With a combination of lattice reduction algorithms we find vectors of norm 4 in all these lattices but in one
pair; only $J_{139} $ and its complex conjugate 
have minimum norm 6.
The quadratic space $(L_K \perp J_{139})^{\#} / (L_K\perp J_{139})  $ is a hyperbolic plane over 
$\F_{47}$. The two isotropic subspaces correspond to isometric even unimodular extremal lattices 
(the isometry is given by $\diag (-I_2,I_{46}) $). 
Let $L$ be one of these lattices. 
By the uniqueness of $L$ it is clear that $L$ is isometric to $P_{48q}$.
To establish an explicit isometry one may consider a vector $\beta \in L_K \leq L$ of norm 12 and
compute the neighbor $L^{(\beta),2}$. This lattice contains a 3-frame, so $L = \Lambda (C)$ for some
ternary extremal code with an automorphism of order 47. 
\eb

\begin{theorem}\label{aut55}
There is no extremal even unimodular lattice $L$ of dimension $48$ such that $\Aut(L) $ contains an
element $\sigma $ of order $55$.
\end{theorem} 

\bew
Then $d=5,p=11$, the type of $\sigma ^5$ is  $11-(4,8)-4$,
 $\det(M) = 11^4$, $\sigma ^5$ acts trivially on $M^{\#}/M$ so $M^{\#}/M \cong \Z[\zeta_{55}]/(1-\zeta _{55}^5)$
and the minimal polynomial for the action of $\sigma $ on $M^{\#}/M$ is $x^4+x^3+x^2+x+1= \Phi_5$.
The automorphism group of the fixed lattice $K$ of $\sigma ^5 $  contains
3 conjugacy classes $a_1,a_2,a_3$
of elements of order 5 with irreducible minimal
polynomial. 
These are the possible candidates for the action of $\sigma $
on $K$. 
They act on the dual quotient $K^{\#}/K$ with minimal polynomial
$$\mu _1:=x^2+8x+1,\ \mu_2 := x^2+4x+1 ,\ \mu_3 = x^4+x^3+x^2+x+1 .$$
So $\sigma $ acts on $M\perp K$ as $\diag (\zeta _{55}, a_3 )$. 
\\
The class group of $F$ is generated by
$\wp = (11, \zeta _{55} +8 ) $, $h_F=10$ and the classes of totally positive units in $\Z_{F^+}$ are represented by
 $U = \{ 1=u_1,u_2,u_3,u_4 \} $.
The ideal
$\overline{\wp} \wp \cap  \Z _{F^+} =  \beta \Z _{F^+}$ has no totally  positive generator
but of course its square is generated by the totally
positive element $\beta ^2$. Also there is some prime element $p_5\in \Z_{F^+}$ dividing 5
such that $p_{55} := \beta p_{5}$ is totally positive.
Therefore the lattice $M$ is isometric to one of the $4\cdot 5$ ideal lattices
$$J_{j,i}:=( \Delta \wp ^{2j+1} , (x,y) \mapsto \trace (u_i \beta ^{-2j} p_{55}^{-1} x \overline{y})) 
 \mbox{ for some } j = 0,\ldots , 4, i = 1,\ldots ,4 .$$
$4\cdot 4$ of these lattices ($j=0,1,3,4$ and all $u_i$) have minimum 6.
\\
 We compute the unimodular overlattices $L$ of
$J_{j,i}\perp L_K$ that are invariant under $\diag (\sigma , a_3 )$.
None of these lattices has minimum 6.
\eb

Similar as in Theorem \ref{aut47} and \ref{aut55} we obtain the following theorem.

\begin{theorem}\label{varphikleiner48} 
Let $L$ be an extremal even unimodular lattice such that $\Aut (L)$ contains some 
element $\sigma $ of order $m$ with $24< \varphi (m) < 48$. 
Then either $m=47$  and $L\cong P_{48q}$, $m=69$ and $L\cong P_{48p}$, $m=120$ and $L\cong P_{48n}$, or $m=132$ and $L\cong P_{48p}$.
\end{theorem}

\bew
The Magma computations 
are similar for all cases and follow the ones described above for the
case $m=47$ and $m=55$. We will always denote by $F:=\Q[\zeta _m]$ the $m$-th
cyclotomic number field, $F^+ := \Q[\zeta _m + \zeta _m^{-1}]$ its maximal real subfield.
The class number of $F^+$ turns out to be 1, the class group of $F$ is cyclic in all cases,
$\wp $ will denote some prime ideal of $\Z _F$  whose class generates the class group of
$F$ and $h_F$ its order, the class number of $F$.
As before, $\Delta $ will denote the inverse different of $\Z_F$, this is always a principal
ideal.
The ideal lattice $M$ will always be of the form $M = L_I (\sigma^d)$ with $dp=m$ for some prime $p$.
In particular we know the $\Z_F$-module 
$M^{\#}/M$ from the computations of the possible
fixed lattices of prime order automorphisms in Section \ref{prime}
as explained above.
The set $U\subset \Z _{F^+}$ will denote some set of units that 
that contains representatives of
all classes of totally positive units.
\\

\noindent
{\bf m=49:} Then $d=p=7$, the  type of $\sigma ^7$ is 
$7-(7,6)-5$, $\det (M) = 7^5$
$\wp=  (197, \zeta _{49} +4 ) $, $h_F=43$,
$\overline{\wp} \wp\cap F^+ = \alpha \Z _{F^+} $ for some totally positive $\alpha $, all totally positive units in  
$\Z_{F^+}$ are squares, so $U=\{ 1 \}$.
A combination of lattice reduction algorithms we find vectors of norm 4 in all 43 lattices
$$J_j:=(\Delta (1-\zeta _{49})^{-1} \wp ^{j} , b_ {\alpha ^{-j}})  \mbox{ for } j = 0,\ldots , 42 .$$
{\bf m=63:}
By Theorem \ref{aut7} there is no automorphism of Type $7-(6,12)-s$, so
$d=21$ and $p=3$ and the type of $\sigma ^{21}$ is
$3-(18,12)-12$ and $M$ is an ideal lattice in the 63-th cyclotomic number field with $\det(M) =3^{12}$. 
The class group of $F$ has order 7 and is generated by any prime ideal dividing 2, e.g. 
$\wp = (2, 1+\zeta_{63}+\zeta_{63}^4+\zeta_{63}^5+\zeta_{63}^6 ) $.
In the maximal real subfield $F^+$ there are unique prime ideals $\wp _7$ and $\wp_3$ dividing $7$ resp. $3$. 
As the class number of $F^+$ is 1, both ideals are principal, $\wp _7$ has some totally positive generator, $\alpha _7$,
but $\wp _3$ doesn't. So the norms of the totally positive elements in $\Z _{F^+}$ that divide 3 are 
powers of $3^{12}$. The lattice 
$(\Z_F, b_{\alpha _7/21} ) $ has determinant $3^{18}$, so 
there is no ideal lattice in the 63rd cyclotomic number field of determinant $3^{12}$.
\\
{\bf m=69:}
Then $d=3$, $p=23$ and $\sigma ^3$ has
 type $23-(2,4)-2$. $\det(M) =23^2$, $h_F = 69$,
$\wp = (139, \zeta _{69} +135 ) $.
$U= \{ 1, u_0 \}$ and
$\overline{\wp} \wp \cap \Z_{F^+}  = \alpha \Z _{F^+}$
for some totally positive element $\alpha \in F^+ $.
So $M$ is isometric to one of the $2\cdot 69$ ideal lattices
$$J_j:=( \Delta \wp ^{j} , (x,y) \mapsto \trace (23 p_3 u \alpha ^{-j}  x \overline{y})) 
 \mbox{ for some } j = 0,\ldots , 68, u \in \{ 1, u_0 \} ,$$
where $p_3$ is some totally positive prime element dividing 3.
$2$ of these lattices ($j=0$ and both $u$)
have minimum 6.
\\
The automorphism group of the fixed lattice $K$ of $\sigma ^{3} $  contains a unique
conjugacy class
of elements of order 3, say represented by $a_1$.
 We compute the unimodular overlattices $L$ of
$M\perp K$ that are invariant under $\diag (\sigma , a_1 )$
and find 24 lattices of minimum 6.
For all these 24 lattices $L$ we compute the
 norm 12-vectors $v$ in the fixed lattice
$K$.
For all lattices there is at least one  $v$
 for which the neighbor $L^{(v),2}$ contains a 3-frame. So $L=\Lambda (C)$ for some
extremal ternary code which turns out to be equivalent to the Pless code $P_{48}$.
\\
{\bf m=75:}
Then $d=15$, $p=5$, the type of $\sigma ^{15}$ is $5-(10,8)-8$, $\det (M) = 5^8$, $h_F = 11$,
$\wp = (151, 146+\zeta _{75}) $, $U=\{ 1,u_0 \}$ and $\overline{\wp} \wp \cap \Z_{F^+}  = \alpha \Z _{F^+}$
Therefore the lattice $M$ is isometric to one of the $2\cdot 11$ ideal lattices
$$J_j:=( \Delta \wp ^{j} , (x,y) \mapsto \trace (75 u \alpha ^{-j} \alpha _{3,5}^{-1} x \overline{y})) 
 \mbox{ for some } j = 0,\ldots , 10, u \in \{ 1, u_0 \} $$
where $\alpha _{35} $ is a totally positive element generating the product of some prime ideal 
over 3 and some prime ideal over 5. 
All 20 non-principal ideal lattices have minimum 6. 
The element $\sigma $ acts on $M^{\#}/M$ as a primitive 15th root of 1.
The automorphism group of the fixed lattice $K=\sqrt{5}E_8$ contains a unique conjugacy class of
automorphisms that act as primitive 15th root of 1, say represented by $a_1$
 We compute the unimodular overlattices $L$ of
$M\perp K$ that are invariant under $\diag (\sigma , a_1 )$.
None of these lattices has minimum 6. 
\\
{\bf m=92:}
Then $d=4$, $p=23$, the type of $\sigma ^{4}$ is $23-(2,4)-2$, $\det (M) = 23^2$, $h_F = 201$,
$\wp = (277, \zeta _{92} +275 ) $.
$U= \{ 1, u_0 \}$ and
$\overline{\wp} \wp \cap \Z_{F^+}  = \alpha \Z _{F^+}$
for some totally positive element $\alpha \in F^+ $.
Therefore the lattice $M$ is isometric to one of the $2\cdot 201$ ideal lattices
$$J_j:=( \Delta \wp ^{j} , (x,y) \mapsto \trace (46 u \alpha ^{-j}  x \overline{y})) 
 \mbox{ for some } j = 0,\ldots , 200, u \in \{ 1, u_0 \} .$$
None of these lattices has minimum 6.
\\
{\bf m=100:}
Then $d=20$, $p=5$, the type of $\sigma ^{20}$ is $5-(10,8)-8$, $\det (M) = 5^8$, $h_F = 5\cdot 11$,
$\wp_5 = (5, 6+\zeta _{100}) $ is an element of order 5 and 
$\wp _7 = (7,1+\zeta _{100}+\zeta _{100}^2-\zeta _{100}^3+\zeta _{100}^4 ) $ an element of order 11 in the class group of $F$. 

Again there are two square classes of totally positive units, so $U=\{ 1,u_0 \}$.
Both ideals 
 $\overline{\wp_5} \wp _5 \cap \Z_{F^+}  = \alpha  _5\Z _{F^+}$
 and $\overline{\wp _7} \wp _7 \cap \Z_{F^+}  = \alpha  _7\Z _{F^+}$ have totally positive generators.
Therefore the lattice $M$ is isometric to one of the $2\cdot 5\cdot 11$ ideal lattices
$$( \Delta \wp _5 ^{j} \wp_7^k , (x,y) \mapsto \trace (50 \alpha _5^{-1} u \alpha _5^{-j} \alpha _7^{-k} x \overline{y})) 
 \mbox{ for some } j = 0,\ldots , 4, k = 0, \ldots , 10,  u \in \{ 1, u_0 \} .$$
All 40 ideal lattices for $j=1,4$ and $k\neq 0$ have minimum 6. 
The element $\sigma $ acts on $M^{\#}/M$ as a primitive 20th root of 1.
The automorphism group of the fixed lattice $K=\sqrt{5}E_8$ contains a unique conjugacy class of
automorphisms that act as primitive 20th root of 1, say represented by $a_1$
 We compute the unimodular overlattices $L$ of
$M\perp K$ that are invariant under $\diag (\sigma , a_1 )$.
None of these lattices has minimum 6. 
\\
{\bf m=108:}
Here $d=36$ and $p=3$ and the type of $\sigma ^{36}$ is
$3-(18,12)-12$ and $M$ is an ideal lattice in the 108-th cyclotomic number field $F$ with $\det(M) =3^{12}$. 
The class group of $F$ has order 19 and is generated by any prime ideal dividing 109, e.g. 
$\wp = (109, 106+\zeta_{108}) $.
In the maximal real subfield $F^+$ there are unique prime ideals $\wp _2$ and $\wp_3$ dividing $2$ resp. $3$. 
As the class number of $F^+$ is 1, both ideals are principal, but none of them has 
a totally positive generator. However $\wp _2\wp _3 = (\alpha _{2,3})$ has some totally positive generator as well 
as $\wp _3^2 = (\alpha _3^2)$. So the norms of the totally positive elements in $\Z _{F^+}$ that divide 3 are 
powers of
$3^{4}$. The lattice 
$(\Z_F, b_{1/18} ) $ has determinant $3^{18}$, so 
there is no ideal lattice in the 108th cyclotomic number field of determinant $3^{12}$.
\\
{\bf m=120:}
The classnumber is $h_F = 4$ and the class group is generated by some prime ideal 
dividing $5$, e.g. $\wp _5 = (5, 2+4\zeta_{120}+\zeta_{120}^2) $. 
In $\Z _F^+$ both prime ideals dividing 3 have totally positive generators, $\alpha _{3a}, \alpha _{3b}$. 
Also the unique prime ideal that divides 2 has a totally positive generator, $\alpha _2$,
whereas the two prime ideals dividing 5 are not generated by totally positive elements, 
their product and their squares are generated by the totally positive elements $\alpha _{5ab}$, $\alpha _{5a}^2$
and $\alpha _{5b}^2$. In particular $\overline{\wp_5} \wp _5\cap F^+$ is not generated by some totally positive element,
so there are only two ideal classes to be considered.
The totally positive units lie in the classes of $1=u_0,u_1,u_2,u_3$. 
\\
Now we have two possibilities for $\det(M)$:  \\
{\bf a)} $d=40$, $p=3$, type of 
 $\sigma ^{40}$ is $3-(16,16)-16$, $\det (M) = 3^{16}$.
Since $\sigma ^{40}$ acts trivially on $K^{\#}/K$ the $\Z_F$-module $M^{\#}/M$ is isomorphic 
to $\Z_F / (1- \zeta _{120}^{40} ) $ and $\sigma $ acts on $M^{\#}/M$ with minimal polynomial $\Phi _{40}$. 
So $M$ is one of the 8 ideal lattices 
$$
( \Delta \wp _5 ^{2j} , (x,y) \mapsto \trace (60 \alpha _5^{-2j} u \alpha _{5ab}^{-1}  x \overline{y})) 
 \mbox{ for some } j = 0,1, u \in \{ 1, u_1, u_2,u_3 \} .$$
Only the lattices for $j=0$ and $u=u_1,u_3$ have minimum 6. 
Both lattices are isometric to the Eisenstein lattice described
by Christine Bachoc \cite{Bachoc}. 
\\
In principle we have two possibilities for the fixed lattice $K$,
 $K\cong \sqrt{3} D_{16}^+$ or $K\cong \sqrt{3} (E_8\perp E_8)$.
The automorphism group of $D_{16}^+$ does not contain elements with minimal polynomial $\Phi _{40}$, so $K\cong \sqrt{3} (E_8\perp E_8)$.
For this lattice $\Aut (E_8\perp E_8)$ contains a unique conjugacy class of such elements,
 represented by, say, $a_1$. 
We now need to compute the $(\sigma , a_1)$ invariant unimodular overlattices of $M\perp K$ for the
two possibilities for $M$. 
A straight forward approach is too memory consuming, 
as we now need to handle invariant submodules 
of $\F _3^{32}$. 
The lattices $L$ are of the form 
$$ L_{\varphi } = \{ (x,y) \in M^{\#} \perp K^{\#} \mid \varphi (x+M) = y+K \} $$
for some isometry  $\varphi : (M^{\#}/M, F_M) \to (K^{\#} /K, -F_K ) $ so that $\sigma \varphi = \varphi a_1 $. 
For any $c\in C_{\Aut (K)}(a_1) = \langle a_1 \rangle $ the lattices $L_{\varphi }$ and $L_{\varphi c}$ are 
isometric. Fixing one such isometry $\varphi _0$, we hence obtain all $L_{\varphi }$ by letting 
$\varphi $ vary in $\{ \varphi _0 u \mid u \in C_{O(K^{\#}/K)} (a_1) / \langle a_1 \rangle \} $. 
To compute the centralizer in the orthogonal group $O(K^{\#}/K)$ we compute the orthogonal elements 
in the endomorphism algebra $\End _{K^{\#}/K}(a_1) \cong \bigoplus _{i=1}^4 \F_{3^4} e_i$. 
As the involution $x\mapsto \overline{x} = (F_K x^{tr} F_K^{-1} ) $ on 
$\End _{K^{\#}/K}(a_1) $ interchanges the idempotents, say, $\overline{e_1} = e_2$, $\overline{e_3}=e_4$, 
the orthogonal elements are of the form $xe_1+\overline{x}^{-1}e_2+ye_3+\overline{y}^{-1} e_4 $,
so $|C_{O(K^{\#}/K)} (a_1) | = 80^2$ and $ C_{O(K^{\#}/K)} (a_1) / \langle a_1 \rangle$ has $160$ coset 
representatives. 
For each lattice $M$ there are exactly $4$ lattices $L_{\varphi }$ that have minimum norm 6.
\\
To obtain an isometry between $L_{\varphi } $ and $P_{48n}$ we 
observe that $\Aut (P_{48n})$ has a unique conjugacy class of automorphism of order 120, say represented by $\tau $. 
We compute a lattice $\Lambda $ such that
$\Lambda /2L_{\varphi }$ is the (4-dimensional) fixed space of the action of $\sigma ^{5}$ on $L_{\varphi }/2L_{\varphi }$.
The dual $D$ of $\Lambda $ (rescaled to be integral) is a lattice of determinant $2^8$, minimum 6 and
kissing number 3264000 and automorphism group of order 5760.
We compute the corresponding sublattice $D_0$  of $P_{48n}$ and an isometry between $D_0$ and $D$.
For all eight lattices $L_{\varphi } $ the isometry maps $L_{\varphi }$ to $P_{48n}$.
\\
{\bf b)} $d=24$, $p=5$, type of 
$\sigma ^{24}$ is $5-(8,16)-8$, $\det (M) = 5^{8}$.
Again $M^{\#}/M \cong \Z_F / (1- \zeta _{120}^{24} ) $ as $\sigma ^{24}$ acts trivially on $K$ and
$\sigma $ acts on $M^{\#} /M$ as a primitive $24$-th root of unity.
So $M$ is one of the 8 ideal lattices 
$$
( \Delta \wp _5 ^{2j} , (x,y) \mapsto \trace (60 \alpha _5^{-2j} u \alpha _{3a}^{-1} \alpha _{3b} ^{-1}  x \overline{y})) 
 \mbox{ for some } j = 0,1, u \in \{ 1, u_1, u_2,u_3 \} .$$
None of these lattices has minimum $\geq 6$.
\\
{\bf m=132:}
Then $d=12$, $p=11$,
$\sigma ^{12}$ has type $11-(4,8)-4$, $\det(M) = 11^4$,
$\wp = (23, \zeta _{132}^2+22\zeta_{132} +8 ) $, $h_F = 11$.
$U= \{ 1, u_0 \}$ and
$\overline{\wp} \wp \cap \Z_{F^+}  = \alpha \Z _{F^+}$
for some totally positive element $\alpha \in F^+ $.
Therefore the lattice $M$ is isometric to one of the $2\cdot 11$ ideal lattices
$$J_j:=( \Delta \wp ^{j} , (x,y) \mapsto \trace (u \alpha ^{-j}  x \overline{y})) 
 \mbox{ for some } j = 0,\ldots , 10, u \in \{ 1, u_0 \} .$$
$10\cdot 2$ of these lattices ($j=1,\ldots,10$ and all $u$)
 have minimum 6.
\\
The automorphism group of the fixed lattice $K$ of $\sigma ^{12} $  contains a unique
conjugacy class
of elements of order 12  with irreducible minimal
polynomial, say represented by $a_1$.
 We compute the unimodular overlattices $L$ of
$M\perp K$ that are invariant under $\diag (\sigma , a_1 )$
and find 240 lattices of minimum 6.
For all these 240 lattices $L$ we compute the 50
orbits of $\langle a_1 \rangle$ on the norm 12-vectors $v$ in the fixed lattice
$L_K$ and find
 one orbit for which the neighbor $L^{(v),2}$ contains a 3-frame. So $L=\Lambda (C)$ for some
extremal ternary code which turns out to be equivalent to the Pless code $P_{48}$.
\eb

\section{The automorphism group of the 72-dimensional extremal lattice} 

In the paper \cite{dim72} I describe the construction of some extremal even unimodular lattice $\Gamma $ 
of dimension 72 such that 
$G:=\Aut(\Gamma )$ contains the subgroup ${\mathcal U} \cong (\SL_2(25) \times \PSL_2(7) ) :2$.
As in Section \ref{aut48} we may prove that ${\mathcal U}$ is indeed the full automorphism 
group of $\Gamma $.
With the same strategy as in the proof of \cite[Theorem 5.3]{cyclo} one shows the following lemma.

\begin{lemma}\label{N25}
${\mathcal U} = N_G(\SL_2(25))$.
\end{lemma}

\bew
For the proof we first note that the natural character of the normal subgroup $N:=\SL_2(25)$ 
is a multiple of one of the two irreducible characters of degree 12. 
The outer automorphism group of $\SL_2(25)$ has order $4$, but only a subgroup of order 2 
stabilises this character. Since ${\mathcal U}$ contains an element inducing this non-trivial 
outer automorphism, it is enough to show that ${\mathcal U}$ contains the centraliser $C:=C_G(N)$.
Let $A$ be the endomorphism algebra of $N$, so $A \cong {\mathcal Q}_{\infty,5}^{3\times 3}$. 
Then $C$ is a subgroup of the unit group $A^*$; moreover $C$ stabilizes $\Gamma $ and the 
${\mathcal U}$-invariant quadratic form. 
Let $v\neq 0$ be any vector in $\Gamma $ and let $L_v := vA \cap \Gamma \leq \Gamma $. 
The dimension of $L_v$ divides $\dim (A) = 36$ and $C$ is a subgroup of $\Aut(L_v)$.
One finds a vector $v$ such that $\dim (L_v) = 36$ and $\Aut(L_v) \cong \pm \PSL_2(7) $.
Therefore $C\cong \pm \PSL_2(7) \leq {\mathcal U}$.
\eb

Using the classification of finite simple groups and in particular the tables in \cite{HM} 
we now conclude as in Section \ref{aut48} that ${\mathcal U} = \Aut (\Gamma )$. 

\begin{theorem}
Let $\Gamma $ be the extremal even unimodular $72$-dimensional lattice constructed in \cite{dim72}.
Then $\Aut(\Gamma ) \cong (\SL_2(25) \times \PSL_2(7) ) :2$.
\end{theorem}

\bew
Explicit generators for the subgroup ${\mathcal U} \cong (\SL_2(25) \times \PSL_2(7) ) :2$ 
of $G:=\Aut(\Gamma )$ have been constructed in \cite{dim72}. 
Remark 2.4 of this paper also shows that ${\mathcal U}$ is absolutely irreducible and that 
all ${\mathcal U}$-invariant lattices are similar to $\Gamma $.
In particular ${\mathcal U}$ does not fix an orthogonally decomposable lattice and therefore 
$G$ is a primitive maximal finite rational matrix group.
In particular the maximal normal $p$-subgroups of $G$ are given in 
the Theorem of Ph. Hall on \cite[p. 357]{Hup}. None of the  possible groups contains $\SL_2(25)$ as 
an automorphism group. 
The tables in \cite{HM} show that there are no quasisimple groups $H\leq \GL_{72} (\Q )$ that 
contain the group ${\mathcal U}$. Therefore one concludes that $\SL_2(25)$ is normal in $G$,
so by Lemma \ref{N25} we obtain $G={\mathcal U}$.
\eb

\end{document}